\documentclass[10pt]{amsart}
\usepackage{amsmath}
\usepackage{amssymb}
\usepackage{amsfonts}
\usepackage{mathrsfs}
\usepackage{graphicx}
\usepackage{epsfig}

\DeclareMathAlphabet{\mathpzc}{OT1}{pzc}{m}{it}

\theoremstyle{plain}

\newtheorem{prop}{Proposition}[section]

\newtheorem{coro}[prop]{Corollary}
\newtheorem{lemm}[prop]{Lemma}
\newtheorem{theoalph}{Theorem}

\theoremstyle{definition}
\newtheorem{defi}[prop]{Definition}

\theoremstyle{remark}
\newtheorem{rema}[prop]{Remark}

\newtheoremstyle{citing}
  {3pt}
  {3pt}
  {\itshape}
  {}
  {\bfseries}
  {.}
  {.5em}
  {\thmnote{#3}}

\theoremstyle{citing}
\newtheorem*{generic}{}

\newcommand{\partn}[1]{{\smallskip \noindent \textbf{#1.}}}

\newcommand\C{\mathbb{C}}
\newcommand\D{\mathbb{D}}

\newcommand\R{\mathbb{R}}

\newcommand\cB{\mathcal{B}}

\newcommand\cF{\mathcal{F}}

\newcommand\cK{\mathcal{K}}

\newcommand\cO{\mathcal{O}}

\newcommand\cR{\mathcal{R}}

\newcommand\cV{\mathcal{V}}
\newcommand\cW{\mathcal{W}}

\newcommand\tA{\widetilde{A}}
\newcommand\tB{\widetilde{B}}

\newcommand\tD{\widetilde{D}}

\newcommand\tK{\widetilde{K}}

\newcommand\tR{\widetilde{R}}

\newcommand\tU{\widetilde{U}}
\newcommand\tV{\widetilde{V}}
\newcommand\tW{\widetilde{W}}

\newcommand\hh{\widehat{h}}

\newcommand\hU{\widehat{U}}
\newcommand\hV{\widehat{V}}
\newcommand\hW{\widehat{W}}

\newcommand\ov{\overline}

\renewcommand{\=}{ : = }

\DeclareMathOperator{\diam}{diam}
\DeclareMathOperator{\dist}{dist}
\DeclareMathOperator{\Area}{Area}

\DeclareMathOperator{\modulus}{mod}
\newcommand{\CV}{\text{CV}}

\newcommand\sC{\mathscr{C}}
\newcommand\hsC{\widehat{\mathscr{C}}}
\newcommand\hA{\widehat{A}}
\newcommand\hH{\widehat{H}}
\newcommand\hK{\widehat{K}}

\newcommand\hdelta{\widehat{\delta}}

\newcommand\hchi{\widehat{\chi}}
\newcommand\hhV{\widehat{\hV}}
\newcommand\hhW{\widehat{\hW}}
\newcommand\hcW{\widehat{\cW}}
\newcommand\hhcW{\widehat{\hcW}}
\newcommand\tp{\widetilde{p}}
\newcommand\tq{\widetilde{q}}
\newcommand\tdelta{\widetilde{\delta}}
\newcommand\teta{\widetilde{\eta}}

\newcommand\tPsi{\widetilde{\Psi}}
\newcommand\tXi{\widetilde{\Xi}}
\newcommand\tchi{\widetilde{\chi}}

\DeclareMathOperator{\Jac}{Jac}

\newcommand{\map}{R}
\newcommand\CC{\ov{\C}}
\newcommand{\Crit}{\text{Crit}}
\newcommand{\CJ}{\Crit(\map) \cap J(\map)}
\newcommand{\CVJ}{\CV(\map) \cap J(\map)}
\newcommand\modul{\mathbf{m}}

\newcommand{\CL}{Connecting Lemma}
\newcommand{\KDT}{Koebe Distortion Theorem}
\newcommand{\BCn}{Backward Contraction}
\newcommand{\BCg}{Backward Contracting}
\newcommand{\BC}[2]{$($#1,#2$)$\nobreakdash-\BCg}
\newcommand{\UPC}{Univalent Pull-back Condition}
\newcommand{\UP}[2]{$($#1,#2$)$\nobreakdash-\UPC}
\newcommand{\deffunr}{$r : (0, \delta_0) \to (1, + \infty)$}
\newcommand\dconst{\mathbf{r}}
\newcommand\clae{conformal Lebsgue a.e.}
\newcommand\cc{connected component}
\newcommand\GL{Gluing Lemma}

\begin{document}
\title[A connecting lemma for rational maps.]{A connecting lemma for rational maps \\ satisfying a no growth condition.}
\author[J. Rivera-Letelier]{Juan Rivera-Letelier$^\dag$}
\thanks{Partially supported by Fundaci\'on Andes and a "Beca Presidente de la Rep\'ublica" of Chile.}
\address{\dag Departamento de Matem\'aticas, Universidad Cat\'olica del Norte, Casilla~1280, Antofagasta, Chile.}
\email{rivera-letelier@ucn.cl}
\date{\today}

\begin{abstract}
We introduce and study a non uniform hyperbolicity condition for complex rational maps, that does not involve a growth condition.
We call this condition \textsf{\textit{\BCn}}.
We show this condition is weaker than the Collet-Eckmann condition, and than the summability condition with exponent~$1$.

Our main result is a connecting lemma for \BCg{} rational maps, roughly saying that we can perturb a rational map to connect each critical orbit in the Julia set with an orbit that does not accumulate on critical points.
The proof of this result is based on Thurston's algorithm and some rigidity properties of quasi-conformal maps.
We also prove that the Lebesgue measure of the Julia set of a \BCg{} rational map is zero, when it is not the whole Riemann sphere.
The basic tool of this article are sets having a Markov property for backward iterates, that are holomorphic analogues of nice intervals in real one-dimensional dynamics.
\end{abstract}

\maketitle
\setcounter{tocdepth}{1}
\tableofcontents

%
%

\section{Introduction.}
We consider rational maps $R : \CC \to \CC$ of degree at least~$2$, viewed as dynamical systems.
We introduce and study a non uniform hyperbolicity condition that we call ``\BCn''.
This condition is defined in~\S\ref{ssi:BC}, and in~\S\ref{ssi:summ and CE} we relate this condition with some well known conditions of non uniform hyperbolicity, such as the Collet-Eckmann and summability conditions.
Our main results are stated in~\S\S\ref{ssi:connecting lemma},~\ref{ssi:expansion away}.
After some notes in~\S\ref{ssi:notes and references}, we describe the organization of this article in~\S\ref{ssi:organization}.

\subsection{\BCn.}\label{ssi:BC}
Let $R$ be a rational map of degree at least~$2$.
We denote by $\Crit(R)$ the set of critical points of $R$ and by $\CV(R) = R(\Crit(R))$ the set of its critical values.
We denote by $J(R)$ the Julia set of $R$; see~\cite{CG,Milbook} for background in holomorphic dynamics.
For $c \in \CJ$ and $\delta > 0$ we denote by $\tB(c, \delta)$ the connected component of $R^{-1}(B(R(c), \delta))$ that contains $c$.

Given $\delta' > \delta > 0$ we will say that $R$ is \textsf{\textit{\BC{$\delta$}{$\delta'$}}}, if for every $c \in \CJ$, every integer $n \ge 1$ and every connected component $W$ of $R^{-n}(\tB(c, \delta'))$, we have that
$$
\dist(W, \CV(R)) \le \delta
\text{ implies }
\diam(W) < \delta.
$$
Given a constant $\dconst > 1$ we will say that $R$ is \textsf{\textit{\BCg{} with constant $\dconst$}}, if for every $\delta > 0$ sufficiently small $R$ is \BC{$\delta$}{$\dconst\delta$}.
Our main results are for rational maps that are \BCg{} with some constant only depending on the rational map.

To better compare the \BCn{} condition with other conditions of non-uniform hyperbolicity, we make the following definition.
Given $\delta_0 > 0$ and a function \deffunr, we will say that $R$ is \textsf{\textit{\BCg{} with function~$r$}}, if for every $\delta \in (0, \delta_0)$ the rational map~$R$ is \BC{$\delta$}{$\delta r(\delta)$}.

\subsection{Collet-Eckmann and summability conditions.}\label{ssi:summ and CE}
Let $R$ be a rational map of degree at least~$2$.
We say that a critical value~$v$ of $R$ is \textsf{\textit{exposed}} if its forward orbit does not contain critical points.

We say that~$R$ satisfies the \textsf{\textit{Collet-Eckmann condition}} if for every exposed critical value~$v$ in $J(R)$, the derivatives $|(R^n)'(v)|$ grow exponentially with~$n$.
The Collet-Eckmann condition was introduced in~\cite{CE83}, in the context of unimodal maps.
It has been extensively studied for complex rational maps, see~\cite{Pr2,Pr1,GSmCE,SmCE,GSw,PRS,Asp, PR} and references therein.

On the other hand, we say that $R$ satisfies the \textsf{\textit{summability condition with exponent}} $\beta > 0$, if for every exposed critical value $v$ of $R$ in $J(R)$, we have
$$
\sum_{n \ge 0} |(R^n)'(v)|^{-\beta}
<
\infty.
$$
The summability condition was introduced in~\cite{NS}, in the context of unimodal maps.
In the complex setting it has been studied in~\cite{Pr1,GSm,PrUparabolic,BvS}.
See also~\cite{Av, Levsumm, Mak}.

In the literature the Collet-Eckmann and summability conditions impose the non existence of parabolic cycles (see however~\cite{PrUparabolic}).
We do not impose this restriction here.

\begin{theoalph}\label{t:summ and CE}
For a rational map $R$ of degree at least~$2$, the following properties hold.
\begin{enumerate}
\item[1.]
If~$R$ satisfies the summability condition with exponent~$1$, then for every $\dconst > 1$ the rational map $R$ is \BCg{} with constant $\dconst$.
\item[2.]
If $R$ satisfies the summability condition with exponent $\beta \in (0,1)$, then there is $\delta_0 > 0$ and a function \deffunr{} satisfying,
$$
\int_0^{\delta_0} (r(\delta))^{- \beta /(1 - \beta)} \frac{d\delta}{\delta}
<
\infty,
$$
such that $R$ is \BCg{} with function~$r$.
\item[3.]
If $R$ satisfies the Collet-Eckmann condition, then there are $\alpha \in (0, 1]$, $C > 0$ and $\delta_0 > 0$ such that $R$ is \BCg{} with function \deffunr, defined by $r(\delta) \= C\delta^{-\alpha}$.
\item[4.]
If the closure of the forward orbit of each critical value of $R$ in $J(R)$ contains no critical point, then there is a constant $C > 0$ such that $R$ is \BCg{} with function \deffunr, defined by $r(\delta) \= C \delta^{-1}$.
\end{enumerate}
\end{theoalph}
In~\cite{RS} it is shown that if $P$ is a real polynomial having only real critical points and such that for every exposed critical value $v$ of $P$ in $J(P)$ we have $|(P^n)'(v)| \to + \infty$ as $n \to \infty$, then for every $\dconst > 1$ the polynomial $P$ is \BCg{} with constant~$\dconst$.

A rational map satisfying the ``Topological Collet-Eckmann condition'' is not in general \BCg.
However, it is shown in~\cite{PR} that these rational maps satisfy a strong form of the \BCn{} property.
\subsection{A connecting lemma.}\label{ssi:connecting lemma}
Given a rational map $R$ of degree at least~$2$ we will denote by $\Crit'(R)$ the set of those critical points of $R$ in $J(R)$, whose image by $R$ is an exposed critical value.

The following is the main result of this article.
\begin{generic}[\CL]
Let $R$ be a rational map of degree at least~$2$ and such that $J(R) \neq \CC$.
Then there is a constant $\dconst > 1$ only depending on $R$, such that if $R$ is \BCg{} with constant~$\dconst$, then for every sufficiently small $\delta > 0$ the following property holds.

For each $c \in \Crit'(R)$ let $v_c \in \CC$ be such that $\dist(v_c, R(c)) < \delta$ and such that for every $n \ge 1$ we have $R^n(v_c) \not \in \tB(\Crit'(R), 2\delta)$.
Let $\xi : \CC \to \CC$ be a quasi-conformal homeomorphism that coincides with the identity outside $\tB(\Crit'(R), \delta)$ and such that for every $c \in \Crit'(R)$ we have $\xi(R(c)) = v_c$.
Then there is a continuous map $h : \CC \to \CC$ and a rational map $Q$ of the same degree as $R$, such that the map $\tR \= \xi \circ R : \CC \to \CC$ satisfies,
$$
h \circ \tR \equiv Q \circ h
\text{ on }
\CC.
$$
Moreover $h$ is close to the identity and $Q$ is close to $R$, as $\delta$ is close to $0$.
\end{generic}
We state a strengthened version of this result as Theorem~\ref{t:strong connecting lemma} in~\S\ref{ss:strong connecting lemma}.
We show in particular that there are many choices for the points $v_c$, so that this result is non vacuous.

To prove this result we show that Thurston's algorithm converges in a very specific situation.
We recall this algorithm in~\S\ref{ss:Thurston's algorithm}, see also~\cite{DHalgorithm,HS}.
In the proof we use a rigidity property that was originally used by J.C.~Yoccoz in his (unpublished) proof that the Mandelbrot is locally connected at non renormalizable parameters.
Similar results have also been used in~\cite{Sester,Khan,Hai,Pil}.
\subsection{Expansion and Lebesgue measure of Julia sets.}\label{ssi:expansion away}
Let $R$ be a rational map of degree at least~$2$.
We show that there is a constant $\dconst > 1$ only depending on~$R$, such that if $R$ is \BCg{} with constant~$\dconst$, then~$R$ has no Siegel disks nor Herman rings.
If moreover $J(R) \neq \CC$, then we show that~$R$ cannot have Cremer cycles, see Corollary~\ref{c:expansion away 2}.

Moreover we prove the following result about Lebesgue measure of Julia sets.
See the survey article~\cite{U} for related results.
\begin{theoalph}\label{t:measure}
Let $R$ be a rational map of degree at least~$2$.
Then there is a constant $\dconst > 1$ only depending on $R$, such that if $R$ is \BCg{} with constant~$\dconst$, then the following properties hold.
\begin{enumerate}
\item[1.]
If $J(R) \neq \CC$, then $J(R)$ has zero Lebesgue measure.
\item[2.]
If $J(R) = \CC$ then there is a set of full Lebesgue measure of points in $\CC$ whose forward orbit accumulates on a critical point of $R$.
\end{enumerate}
\end{theoalph}
In~\cite{ims} it is also shown that, if $R$ is a polynomial rational map as in the theorem, then every connected component of the Julia set of $R$ is locally connected.
\subsection{Notes and references.}\label{ssi:notes and references}
This article is part of~\cite{ims}, in revised form.
A rational map satisfies the Decay of Geometry condition in the sense of~\cite{ims}, if for every $\dconst > 1$ it is \BCg{} with constant~$\dconst$.

The proof of Theorem~\ref{t:summ and CE} is based on a technique to control distortion along backward trajectories, called \textsf{\textit{shrinking neighborhoods}}, that was introduced by F.~Przytycki in~\cite{Pr2}.

The converse of part~$4$ of Theorem~\ref{t:summ and CE} is easily seen to hold.

The proof of the results stated in~\S\ref{ssi:expansion away} are based on a strengthened version of the \CL, that we state as Theorem~\ref{t:strong connecting lemma} in~\S\ref{ss:strong connecting lemma}.
This result is only used to prove some expansion properties of forward orbits that do not accumulate on critical points, see~\S\ref{ss:expansion away}.
It would be very desirable to have a more direct proof of these facts.
Furthermore, is not clear to me if a rational map that is \BCg{} with constant~$\dconst$, for every $\dconst > 1$, can have a Cremer periodic point.

\subsection{Organization}\label{ssi:organization}
We briefly describe the organization of this article.

After some preliminaries in~\S\ref{s:preliminaries}, we introduce ``nice sets'' and ``nice nests'' in~\S\ref{s:nice sets and nests}.
These concepts are basic for the rest of this article.
Nice sets are holomorphic analogues of ``nice intervals'' in real one dimensional dynamics.

In \S\S\ref{s:Rigidity}, \ref{s:Thurston's algorithm} we give a sufficient condition under which Thurston's algorithm converges in a very specific situation.

In \S\ref{s:selfsimilar} we give an equivalent formulation of the \BCn{} condition and we show how to construct nice sets and nests for rational maps satisfying this condition.

In \S\ref{s:connecting lemma} we state and prove a strengthened version of the \CL{}.
The proof is based on some key area estimates and on the results on Thurston's algorithm.
In~\S\ref{s:global} we derive some consequences of these results, and in particular we prove Theorem~\ref{t:measure}.

Appendix~\ref{a:summ and CE} contains the proof of Theorem~\ref{t:summ and CE} and in Appendix~\ref{a:qc} we collect some basic facts about quasi-conformal maps.
\subsection{Acknowledgments.}
I am grateful to J.C.~Yoccoz for stimulating conversations, various suggestions and for the careful reading of an early version of this paper.
I am also grateful to J.~Kiwi for sharing insights and for useful conversations, suggestions and comments.
I would like to express my gratitude to P.~Roesh, for useful conversations and comments, and to Universit\'e de Lille~$1$.
Finally, I would like to thank to W.~Shen for a very useful comment.

Figure~\ref{f:McMullen} was done with one of C.~McMullen's programs.

\section{Preliminaries.}\label{s:preliminaries}
For~2 numbers $A$ and $B$, $A \sim B$ and $A = {\mathcal O}(B)$ means $C^{-1}B < A < CB$ and $A < CB$ for some implicit constant $C > 0$, respectively.

We endow the Riemann sphere $\CC$ with the spherical metric.
Distances, diameters, balls and derivatives are all taken with respect to the spherical metric. 
For $z \in \CC$ and $\delta > 0$ we denote by $B(z, \delta)$ the open ball centered at~$z$ and of radius $\delta$.

\subsection{Critical points.}
Given a rational map $R$, we denote by $\Crit(R) \subset \CC$ the set of its critical points and by $\CV(R) = R(\Crit(R))$ the set of its critical values.
When there is no danger of confusion we denote $\Crit(R)$ and $\CV(R)$ just by $\Crit$ and $\CV$, respectively.
For a critical point $c \in \Crit(R)$ we denote by $\mu_c$ the local degree of $R$ at~$c$ and put $\mu_{\max} \= \max \{ \mu_c \mid c \in \CJ \}$.

To treat the case when a critical point in $J(R)$ is mapped to another critical point by some iterate of~$R$, we consider a block $c_1, \ldots, c_k$ of critical points, such that $c_j$ mapped to $c_{j + 1}$ by some iterate of $R$ and maximal with these properties, as a single critical point of multiplicity $\mu_{c_1} \cdot ... \cdot \mu_{c_k}$.
For example the critical value of this sequence means $R(c_k)$.
With this convention we assume that no critical point in $J(R)$ is eventually mapped into some other critical point.
In particular we assume that $\CVJ$ is disjoint from $\Crit(R)$.

Recall that for a critical point $c \in J(R)$ and $\delta > 0$ small, the set  $\tB(c, \delta)$ is the connected component of $R^{-1}(B(R(c), \delta))$ that contains~$c$.
Thus
$$
R(\tB(c, \delta)) = B(R(c), \delta)
\text{ and }
\diam(\tB(c, \delta)) \sim \delta^{\frac{1}{\mu_c}}.
$$

\subsection{Pull-backs.}
Given a subset $V$ of $\CC$ and a non negative integer $n$, each connected component of $R^{-n}(V)$ is called \textsf{\textit{a pull-back of $V$ by $R^n$}}.
Note that here $V$ is not assumed to be connected.
Moreover a subset $W$ of $\CC$ is called a \textsf{\textit{pull-back of}} $V$ if there is a non negative integer~$n$ such that $W$ is a pull-back of $V$ by $R^n$.

It is easy to see that for every rational map~$R$ that is \BCg{} with some constant~$\dconst > 1$, the set $\CV(R)$ is disjoint from the parabolic periodic points of~$R$.
Hence, by the Fatou-Sullivan classification of connected components of the Fatou set~\cite{CG,Milbook}, there is a neighborhood of $\CVJ$ that is disjoint from the forward orbits of critical points not in $J(R)$.
We will implicitly assume that all neighborhoods of points in $\CVJ$ or $\CJ$ that we consider, are sufficiently small to be disjoint from the forward orbits of critical points not in~$J(R)$.

\subsection{Distortion.}\label{ss:distortion}
The following is a version of \KDT{} in $\CC$, endowed with the spherical metric.
It can be easily deduced from the usual \KDT{}, see for example~\cite{Pom}.
\begin{generic}[\KDT]
For each $\varepsilon \in (0, \tfrac{1}{2}]$ there is a constant $D > 1$ such that the following property holds.
For every $r \in (0, \diam(\CC))$, every $z_0 \in \CC$ and every conformal map $\varphi : B(z_0, r) \to \CC$ such that
$$
\diam(\CC \setminus \varphi(B(z_0, r))) > \tfrac{1}{2} \diam(\CC),
$$
the distortion of $\varphi$ on $B(z_0, \varepsilon r)$ is bounded by $D$.
Moreover $D = 1 + \cO(\varepsilon)$.
\end{generic}
We will use this result in the following situation.
From now on, given a rational map $R$ of degree at least~$2$ we fix a coordinate on $\CC$ for which $\diam(J(R)) > \tfrac{1}{2} \diam(\CC)$, so that the hypothesis of the lemma below is satisfied.
So, if $r_K > 0$ is as in the lemma below, then for every $r \in (0, r_K)$, every $z_0 \in \CC$, every positive integer $n$ and every pull-back $W$ of $B(z_0, r)$ by $R^n$, we have $\diam(\CC \setminus W) > \tfrac{1}{2} \diam(\CC)$.
Therefore, if $R^n$ is univalent on $W$, then we can apply \KDT{} to $\varphi = R^n|_W^{-1} : B(z_0, r) \to W$.
\begin{lemm}\label{l:backward size}
Let $R$ be a rational map of degree at least~$2$ such that $\diam(J(R)) > \tfrac{1}{2} \diam(\CC)$.
Then there is $r_K \in (0, \diam(\CC))$ such that for every $r \in (0, r_K)$, every $z_0 \in \CC$ and every pull-back $W$ of $B(z_0, r)$ by $R$, we have
$$
\diam(\CC \setminus W) > \tfrac{1}{2} \diam(\CC).
$$
\end{lemm}
\begin{proof}
As the repelling periodic points of $R$ are dense in $J(R)$ \cite{CG,Milbook}, we can find disjoint forward invariant sets $\cO$ and $\cO'$, consisting of finitely many periodic points of $R$ and such that
\begin{equation}\label{e:backward size}
\diam(\cO), \diam(\cO') > \tfrac{1}{2} \diam(\CC).
\end{equation}
Let $r_K > 0$ be less than the distance between $\cO$ and $\cO'$.
Then for every $r \in (0, r_K)$ and $z_0 \in \CC$, the ball $B(z_0, r)$ is disjoint from either $\cO$ or $\cO'$.
As $\cO$ and $\cO'$ are forward invariant by $R$, it follows that every pull-back $W$ of $B(z_0, r)$ by $R$ is disjoint from either $\cO$ or $\cO'$.
Then the assertion of the lemma follows from~\eqref{e:backward size}. 
\end{proof}

\section{Nice sets and nests.}\label{s:nice sets and nests}
In this section we introduce ``nice sets'' and ``nice nests'', in~\S\ref{ss:nice sets} and~\S\ref{ss:nice nests} respectively.
We first consider some preliminary remarks in~\S\ref{ss:maximal invariant}.
We end this section by describing a general way to construct nice sets and nests in~\S\ref{ss:criterion nice}.

We fix throughout all this section a rational map $R$ of degree at least~2.
\subsection{Maximal invariant sets.}\label{ss:maximal invariant}
We will assume that there is at least one critical point of~$R$ in~$J(R)$.
\begin{defi}
Let $V$ be a neighborhood of $\CJ$ such that every connected component of $V$ contains exactly one critical point in $J(R)$.
Then we define,
$$
K(V) \= \{ z \in \CC \mid \text{for every integer $j \ge 0$ we have } R^j(z) \not \in V \}.
$$
\end{defi}
Let $V$ be as in the definition.
Then the set $K(V)$ is compact, forward invariant by~$R$ and disjoint from~$V$.
On the other hand, Montel's Theorem implies that the interior of $K(V)$ is contained in $\CC \setminus J(R)$; see~\cite{CG,Milbook}.
If $V'$ is a neighborhood of $\CJ$ such that each connected component contains a unique element of $\CJ$, then $V' \subset V$ implies $K(V) \subset K(V')$.

If $W$ is a connected component of $\CC \setminus K(V)$ not intersecting $\CJ$, then $R(W)$ is also a connected component of $\CC \setminus K(V)$ and the map $R : W \to R(W)$ proper.
When $V$ is disjoint from the forward orbits of critical points not in $\CJ$, it follows that $W$ can only contain critical points of $R$ in $J(R)$.
Therefore, for each connected component $W$ of $\CC \setminus K(V)$ there is an integer $m_W \ge 0$ and $c(W) \in \CJ$ so that $R^{m_W}$ maps $W$ biholomorphically onto the connected component of $\CC \setminus K(V)$ containing $c(W)$.
\subsection{Nice sets.}\label{ss:nice sets}
\begin{defi}
For $c \in \CJ$ consider a simply-connected neighborhood $V^c$ of $c$, disjoint form the forward orbits of critical points not in $J(R)$.
Moreover suppose that the closures of the sets $V^c$ are pairwise disjoint and put $V \= \cup_{c \in \CJ} V^c$.
Then we say that $V$ is a \textsf{\textit{nice set for}} $R$ if for every integer $n \ge 1$ and every connected component $W$ of $R^{-n}(V)$, we have either
$$
\ov{W} \cap \ov{V} = \emptyset
\ \text{ or } \
\ov{W} \subset V.
$$
\end{defi}
By definition of nice sets, each pull-back of a nice $V$ of $R$ is disjoint form the critical points of $R$ not in $J(R)$. 
Note moreover that for distinct pull-backs $W$ and $W'$ of a nice set $V$ of $R$, we have either
$$
\ov{W} \cap \ov{W'} = \emptyset,
\
\ov{W} \subset W'
\ \text{ or } \
\ov{W'} \subset W.
$$

M.~Martens introduced the concept of \textsf{\textit{nice interval}} in~\cite{Mar}, in the context of interval maps.
An interval is said to be nice if the forward orbit of every point in its boundary is disjoint from the interval itself.
The following lemma shows that nice sets satisfy an analogous property.
\begin{lemm}\label{l:nice sets forward}
Let $V = \cup_{c \in \CJ} V^c$ be a nice set of $R$.
Then for every point $z \in \partial V$ and every integer $n \ge 1$ we have $R^n(z) \not \in \ov{V}$.
So we have the following properties.
\begin{enumerate}
\item[1.]
For each critical point $c \in \CJ$, the set $V^c$ is equal to the connected component of $\CC \setminus K(V)$ that contains $c$.
In particular $\partial V^c \subset K(V)$.
\item[2.]
For every connected component $W$ of $\CC \setminus K(V)$ there is an integer $m_W \ge 0$ and a critical point $c(W) \in \CJ$ such that 
$$
R^{m_W} : W \to V^{c(W)}
$$
is a biholomorphism.
In particular $W$ is simply-connected and therefore $K(V)$ is connected.
\end{enumerate}
\end{lemm}
\begin{proof}
Suppose by contradiction that for some $z \in \partial V$,  $n \ge 1$ and $c \in \CJ$ we have $R^n(z) \in \ov{V^c}$.
Then the pull-back $W$ of $V^c$ by $R^n$ to $z$ is such that $z \in \ov{W} \cap \partial{V}$. 
Thus $\ov{W} \cap \ov{V} \neq \emptyset$ and $\ov{W} \not\subset V$, which contradicts the hypothesis that $V$ is a nice set.

As remarked above, if $W$ is a connected component of $\CC \setminus K(V)$ different from $V^c$, for each $c \in \CJ$, then $R(W)$ is also a connected component of $\CC \setminus K(V)$ and the map $R : W \to R(W)$ is proper.
It follows that $R : W \to R(W)$ is a biholomorphism.
The rest of the assertions are easily deduced from this property.
\end{proof}
\subsection{Nice nests.}\label{ss:nice nests}
Given a positive integer $\ell$ we will say that a sequence $(V_0, \ldots, V_\ell)$ of nice sets is a \textsf{\textit{nice nest}} for $R$, if $\ov{V_\ell} \subset V_{\ell - 1}, \ldots, \ov{V_1} \subset V_0$ and if for every $j \in \{0, \ldots, \ell \}$ and every pull-back $W$ of $V_0$ by $R$, we have either
\begin{equation}\label{e:nest condition}
\ov{W} \cap \ov{V_j} = \emptyset
\text{ or }
\ov{W} \subset V_j.
\end{equation}
When $\ell = 1$ (resp. $\ell = 2$), the sequence $(V_0, V_1)$ (resp. $(V_0, V_1, V_2)$) is also called \textsf{\textit{nice couple}} (resp. \textsf{\textit{nice triple}}) \textsf{\textit{for}} $R$. 

If $\ell$ is a positive integer and $(V_0, \ldots, V_\ell)$ is a nice nest for $R$, then for every $j, k \in \{0, \ldots, \ell \}$ and every pull-back $W$ of $V_k$ we have~\eqref{e:nest condition}.
On the other hand, if $V$ is a nice set such that $\ov{V} \subset V_\ell$ and such that $(V_0, V)$ is a nice couple, then $(V_0, \ldots, V_\ell, V)$ is a nice nest.

Let $\ell$ be a non negative integer and let $V_0$ be a nice set if $\ell = 0$, and let $(V_0, \ldots, V_\ell)$ be a nice nest if $\ell \ge 1$.
Given $\modul > 0$, we will say that $(V_0, \ldots, V_\ell)$ is $\modul$-\textsf{\textit{separated}} if there is a nice set $\hV$ containing $\ov{V_0}$ such that,
$$
\min_{c \in \CJ} \modulus(\hV^c \setminus \ov{V_0^c}) \ge \modul,
$$
and such that the sequence $(\hV, V_0, \ldots, V_\ell)$ is a nice nest.
Suppose that $(V_0, \ldots, V_\ell)$ is $\modul$-separated.
Let $j \in \{0, \ldots, \ell\}$, let $W$ be a \cc{} of $\CC \setminus R^{-1}(K(V_j))$ disjoint from $V_\ell$ and let $W_0$ be the pull-back of $V_0$ by $R^{m_W}$ containing $W$.
Then $R^{m_W}$ is univalent on $W_0$ and there is a constant $D(\modul) > 1$ only depending on $\modul$ such that the distortion of $R^{m_W}$ on $W_0$ is bounded by $D(\modul)$.
We have $D(\modul) \to 1$ as $\modul \to + \infty$.
\subsection{A way to construct nice sets and nests.}\label{ss:criterion nice}
The following lemma will be used in~\S\ref{ss:construction nice sets and nests} to produce nice sets and nests.
\begin{lemm}\label{l:nice sets with closures}
Let $\tV_0$ be a neighborhood of $\CJ$ in $\CC$ so that each connected component of $\tV_0$ contains precisely one element of $\CJ$.
Let $\ell$ be an non negative integer and for each $c \in \CJ$ let $V_0^c, \ldots, V_\ell^c$ be simply-connected neighborhoods of $c$ in such a way that
$$
\ov{V_{\ell}^c} \subset V_{\ell -1}^c, \ldots, \ov{V_1^c} \subset V_0^c, \ov{V_0^c} \subset \tV_0, 
$$
and that $V_0^c$ is disjoint from the forward orbits of critical points not in $J(R)$.
If for every $j = 0, \ldots, \ell$ we have $R(\partial V_j^c) \subset K(\tV_0)$, then each of the sets
$$
V_j \= \cup_{c \in \CJ}V_j^c
$$
is a nice set and when $\ell \ge 1$ the sequence $(V_0, \ldots, V_\ell)$ is a nice nest.
\end{lemm}
\begin{proof}
Fix $j \in \{0, \ldots, \ell\}$.
Given $c \in \CJ$ let $W_0 = V_0^c, W_1, \ldots$ be successive pull-backs by $R$ and for $n \ge 1$ let $\tW_n$ be the connected component of $\CC \setminus K(\tV_0)$ that contains $W_n$.
Since $\partial V_j \subset R^{-1}(K(\tV_0))$ it follows that for each $n \ge 1$ we have either $\tW_n \subset V_j$ or $\tW_n \cap V_j = \emptyset$.
So to prove the lemma it is enough to prove that for each $n \ge 0$ we have $\ov{W_n} \subset \tW_n$.

We proceed by induction.
For $n = 0$ just note that $\ov{V_0}$ is contained in $\tV_0$ by hypothesis, so $\ov{W_0} \subset \tW_0$.
Suppose by induction hypothesis that for some integer $n \ge 1$ we have $\ov{W_n} \subset \tW_n$.
If $\tW_n$ does not intersect $\CJ$ then the map $R : \tW_n \to \tW_{n - 1}$ is proper, so we have $\ov{W_n} \subset \tW_n$ in this case by the induction hypothesis.
If $\tW_n$ intersects $\CJ$ then let $\tW_n'$ be the connected component of $R^{-1}(\tW_{n - 1})$ that contains $W_n$, so that $\tW_n' \subset \tW_n$.
By the induction hypothesis, we have $\ov{W_n} \subset \tW_n' \subset \tW_n$.
This completes the induction step and ends the proof of the lemma.
\end{proof}

\section{Rigidity.}\label{s:Rigidity}
After some general definitions and facts in~\S\S\ref{ss:rigid pairs}, \ref{ss:rigid nice couples}, we state the main result of this section in~\S\ref{ss:Rigidity}.
The proof of this result in~\S\ref{ss:proof of Rigidity}.
\subsection{Rigid pairs.}\label{ss:rigid pairs}
Let $U$ be an open and connected subset of $\CC$ whose complement contains at least~$3$ points, and let $E$ be a measurable subset of $U$.
Given a constant $K > 1$ we will say that the pair $(U, E)$ is $K$-\textsf{\textit{rigid}} if there is a constant $C > 0$ such that for every qc homeomorphism $\chi : U \to \chi(U) \subset \CC$ that is \clae{} on $E$ there is a $K$-qc homeomorphism $\hchi : U \to \chi(U)$ whose distance to $\chi$ is at most $C$.
We will also say that the pair $(U, E)$ is $K$-\textsf{\textit{rigid with constant}} $C$.
Here the distance is taken with respect to the hyperbolic metric of $\chi(U)$.
We will say that a pair $(U, E)$ is \textsf{\textit{rigid}} if there is $K > 1$ for which it is $K$-rigid.
Note that if the pair $(U, E)$ is $K$-rigid and if $\varphi : U \to \CC$ is a biholomorphism onto its image, then the pair $(\varphi(U), \varphi(E))$ is $K$-rigid with the same constant.

\begin{rema}\label{r:rigid pairs}
The property that the hyperbolic distance between $\chi$ and $\hchi$ is finite, implies that $\chi$ extends continuously to the closure of $U$ if and only if $\hchi$ does.
In this case the extensions of $\chi$ and $\hchi$ coincide on the boundary of $U$.
This property will be important to apply the \GL{} (see Appendix~\ref{a:qc}.)
\end{rema}
\subsection{Rigid nice couples.}\label{ss:rigid nice couples}
Fix a rational map $R$ of degree at least~$2$.
Then we will say that a nice couple $(\hV, V)$ for $R$ is $K$-\textsf{\textit{rigid with constant}} $C$ (resp. $K$-\textsf{\textit{rigid}}, \textsf{\textit{rigid}}) if for every $c \in \CJ$ the pair $(\hV^c, \hV^c \cap K(V))$ is.
\begin{prop}\label{p:rigid nice couples}
Let $\ell \ge 1$ be an integer and let $(V_0, \dots, V_\ell)$ be a nice nest.
Let $K > 1$ and $C > 0$ be given and suppose that there are $j, k \in \{ 0, \ldots, \ell \}$ such that $j < k$ and such that the nice couple $(V_k, V_j)$ is $K$-rigid with constant $C$.
Then for each $j' \in \{0, \ldots, j\}$ the following property holds.
\begin{enumerate}
\item[1.]
For every $k' \in \{k, \ldots, \ell \}$ the nice couple $(V_{k'}, V_{j'})$ is $K$-rigid with constant~$C$.
\item[2.]
For every $k' \in \{0, \ldots, \ell \}$ such that $k' > j'$ and every $c \in \CJ$ the pair $(V_{k'}^c \setminus \ov{V_{j'}^c}, V_{k'}^c \cap K(V_{j'}))$ is $K$-rigid with constant $C$.
Moreover there is a constant $K' > 1$ only depending on
$$
\modul = \min_{c \in \CJ} \modulus(V_{k'}^c \setminus \ov{V_{j'}^c})
$$
and a constant $C' > 0$ only depending on $K$ and $\modul$, such that the nice couple $(V_{k'}, V_{j'})$ is $KK'$-rigid with constant $C + C'$.
\end{enumerate}
\end{prop}
The proof of this proposition is below.
It depends on the following lemmas.
\begin{lemm}\label{l:smoothing2}
Given $r \in (0,1)$ there is $K'(r) > 1$ such that for every
homeomorphism $\chi : \D \to \D$ that is conformal on $\{ r < |z| <1 \}$, there is a $K'(r)$-qc homeomorphism $\tchi : \D \to \D$ that coincides with $\chi$ on $\{ r^{\frac{1}{2}} < |z| < 1 \}$.
Moreover $K'(r) \to 1$ when $r \to 0$.
\end{lemm}
\begin{proof}
Put $h = - \frac{1}{4\pi} \ln r$ and
\begin{eqnarray*}
E : \{ z \in \C \mid \Im z < h \}  & \to & \D \setminus \{ 0 \} \\
z & \mapsto & r^{\frac{1}{2}} \exp(-2\pi i z).
\end{eqnarray*}
Note that $E$ is a $1$-periodic universal covering map and that $E(\R) = \{ |z| = r^\frac{1}{2} \}$ .
On the other hand, let $E_0 : \{ \Im z < 0 \} \to \chi (\{ |z| < r^\frac{1}{2} \})$ be a $1$-periodic universal covering map and let $\tchi : \{ \Im z < 0 \} \to \{ \Im z < 0 \}$ be a lift of $\chi|_{\{ 0 < |z| < r^\frac{1}{2} \}}$.
Since $\chi$ is conformal on $\{ r < |z| < 1 \}$, it follows that $\tchi$ is conformal on $\{ - h < \Im z < 0 \}$.
By Schwarz reflexion principle it follows that $\tchi$ extends to a homeomorphism defined on $\{ \Im z < h \}$, that is conformal on $\{ | \Im z | < h \}$.

Denote by $f : \R \to \R$ the restriction of $\tchi$ to $\R$, which is a $1$-periodic diffeomorphism.
As $f$ has a conformal extension to $\{ |\Im z| < h \}$, it follows by \KDT{} that there is a constant $K'(r) > 1$, only depending on $r = \exp( - 4\pi h)$, such that
$$
\sup_{x \in \R} \max \{ f'(x), (f'(x))^{-1} \} \le K'(r).
$$
Moreover $K'(r) \to 1$ as $r \to 0$.

Let $\tchi_0 : \{ \Im z < h \} \to \tchi(\{ \Im z < h \})$ be the homeomorphism that coincides with $\tchi$ on $\{ 0 \le \Im z < h \}$ and such that $\tchi_0(x + i y) = f(x) + i y$, for $x \in \R$ and $y \le 0$.
Is easy to see that $\tchi_0$ is of class $C^1$ on $\{ \Im z < 0 \}$ and that the dilatation of $\tchi_0$ at $x + i y$, for $y < 0$, is equal to $\max \{ f'(x), (f'(x))^{-1} \} \le K'(r)$.

As $\tchi_0$ is $1$-periodic, it is a lift of a homeomorphism $\tchi : \D \setminus \{ 0 \} \to \chi(\D \setminus \{ 0 \})$ that coincides with $\chi$ on $\{ r^\frac{1}{2} \le |z| < 1 \}$ and that is $K'(r)$-qc on $\{ 0 < |z| < r^\frac{1}{2} \}$.
As the set $\{ |z| = r^\frac{1}{2} \} \cup \{ 0 \}$ is qc removable, it follows that $\tchi$ extends to a $K'(r)$-qc homeomorphism from $\D$ to $\chi(\D)$.
\end{proof}
\begin{lemm}\label{l:smoothing3}
For each $\modul > 0$ there is a constant $K(\modul) > 1$ such that the following property holds.
Let $U \subset \CC$ be biholomorphic to $\D$ and let $\cK \subset U$ be a compact set such that $A \= U \setminus \cK$ is an annulus whose modulus is at least $\modul$.
Then, for every constant $K_0 > 1$ there is a constant $C(K_0, \modul) > 0$ such that for every homeomorphism $\chi : U \to \chi(U) \subset \CC$ that is $K_0$-qc on $A$, there is a $K_0K(\modul)$-qc homeomorphism $\hchi : U \to \chi(U)$ whose hyperbolic distance to $\chi$ is at most $C(K_0, \modul)$.
Moreover $K(\modul) \to 1$ as $\modul \to \infty$ and for a fixed $K_0 > 1$ we have $C(K_0, \modul) \to 0$ as $\modul \to \infty$.
\end{lemm}
\begin{proof}
Without loss of generality assume that $U = \D$ and that $0 \in \cK$.
Then there is $r \in (0, 1)$ only depending on $\modul$ such that $\{ r < |z| < 1 \}$ is contained in~$A$.
Clearly $r \to 0$ as $\modul \to \infty$.
Let $\chi : \D \to \chi(\D)$ be a homeomorphism that is $K_0$-qc on~$A$.
Then note that the hyperbolic diameter of $\chi(\{ |z| \le r^{\frac{1}{2}} \})$ in $\chi(\D)$ is bounded by a constant $C'(K_0, r)$ only depending on $K_0$ and $r$.
Clearly for a fixed $K_0$ we have $C'(K_0, r) \to 0$ as $r \to 0$.

Let $\psi : \chi(U) \to \D$ be a $K_0$-qc homeomorphism so that $\psi \circ \chi$ is conformal on $A$ and let $\tchi : \D \to \D$ be the $K'(r)$-qc homeomorphism given by Lemma~\ref{l:smoothing2}, for $\chi$ replaced by $\psi \circ \chi$.
Then $\hchi \= \psi^{-1} \circ \tchi$ is a $K_0K'(r)$-qc homeomorphism that coincides with $\chi$ on $\{ r^{\frac{1}{2}} < |z| < 1 \}$.
This proves the lemma with constants $K(\modul) \= K'(r)$ and $C(K_0, \modul) \= C'(K_0, r)$.
\end{proof}
\begin{proof}[Proof of Proposition~\ref{p:rigid nice couples}]
\

\partn{1}
As the set $V_{k'}^c \cap K(V_{j'}^c)$ contains $V_{k'}^c \cap K(V_j^c)$, is enough to prove the assertion when $j' = j$ and $k' > k$.
Note that the connected components of $V_{k'}^c \setminus K(V_k)$ cover $V_{k'}^c \setminus K(V_j)$.
Given such a \cc{} $W$ denote by $\phi_W$ the inverse of $R^{m_W}|_{V_{k}^{c(W)}}$.
Then $\chi \circ \phi_W : V_k^{c(w)} \to \chi(W)$ is a qc homeomorphism that is \clae{} on $K(V_j)$.
As the nice couple $(V_k, V_j)$ is $K$-rigid with constant $C$, there is a $K$-qc homeomorphism $\psi_W : W \to \chi(W)$ whose hyperbolic distance to $\chi|_W$ is at most $C$.
Then the assertion follows by the Gluing Lemma.

\partn{2}
That the pair $(V_{k'}^c \setminus \ov{V_{j'}^c}, V_{k'}^c \cap K(V_{j'}))$ is $K$-rigid with constant $C$, can be proved in the same way as in part~$1$.
The second assertion is then a direct consequence of Lemma~\ref{l:smoothing3}.
\end{proof}
\subsection{Rigidity.}\label{ss:Rigidity}
For a subset $X$ of $\C$ we denote by $\Area(X)$ the area of $X$ with respect to the Euclidean metric of $\C$.
Given a subset $U$ of $\CC$ biholomorphic to $\D$ and a subset $E$ of $U$, put
\begin{equation}\label{e:density}
\| E \|_U \= \sup_{\varphi} \Area(\varphi(E)),
\end{equation}
where the supremum is taken over all biholomorphisms $\varphi: U \to \D$.

The following is the main result of this section.
\begin{generic}[Rigidity]
Let $R$ be a rational map of degree at least~$2$ and let $K > 1$ be given.
Then there are constants $\modul > 0$ and $\varepsilon > 0$ such that if $(\hhV, \hV, V)$ is a nice triple for $R$, satisfying 
$$
\min_{c \in \CJ} \modulus(\hV^c \setminus \ov{V^c}) \ge \modul
\text{ and }
\left\| \hhV^c \setminus K(\hV) \right\|_{\hhV^c} < \varepsilon,
$$
then the nice couple $(\hV, V)$ is $K$-rigid.
\end{generic}
The proof of this statement is in \S\ref{ss:proof of Rigidity}.
The following general lemma describes the way in which~\eqref{e:density} is used.
\begin{lemm}\label{l:smoothing1}
Given constants $C_0 > 0$ and $K_0$, $K_1 > 1$, there is $\varepsilon > 0$ such that the following property holds.
Let $U$, $U' \subset \CC$ be biholomorphic to $\D$ and let $E \subset U$ be such that $\| E \|_{U} < \varepsilon$.
Then for every $K_1$-qc homeomorphism $\chi : U \to U'$ that is \clae{} outside $E$, there is a $K_0$-qc homeomorphism $\hchi : U \to U'$ whose hyperbolic distance to $\chi$ is bounded by~$C_0$.
\end{lemm}
\begin{proof}
Without loss of generality we assume $U' = U = \D$.
Consider a tiling of $\D$ by regular hyperbolic hexagons $( H_j )_{j \ge 0}$ with straight angles, as in figure~\ref{f:McMullen}.
\begin{figure}[htb]
\begin{center}
\psfig{file=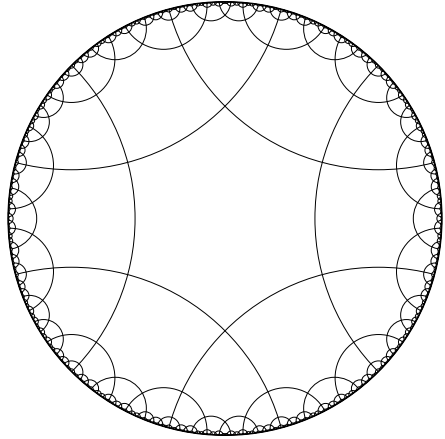, width = 2.5in}
\caption{The vertices of the hexagons from a uniformly distributed set in $\D$.}
\label{f:McMullen}
\end{center}
\end{figure}

Let $H$ by a hyperbolic hexagon with right angles that is centered at $z = 0$ and denote by $\cV$ the set of its vertices (recall that all such hexagons are isometric).
Let $\delta > 0$ be small enough so that for every map $\chi^0 : \cV \to \D$ such that $\max \{ |\chi^0(v) - v| \mid v \in \cV \} \le \delta$, the following properties hold.

$\bullet$
The geodesic segments joining the images by $\chi^0$ of consecutive vertex of $H$, form a hexagon $\hH$.

$\bullet$
Consider the extension $\chi^1 : \partial H \to \hH$ of $\chi^0$, determined by the property that the distance between points on the same side of $\partial H$ is multiplied by a constant factor, that only depends on the side of $\partial H$.
Then $\chi^1$ extends to a homeomorphism from $\ov{H}$ to $\ov{\hH}$, that is $K_0$-qc on the interior of $H$ and whose hyperbolic distance to the identity is at most $C_0 / 2$.

\

Moreover we assume that $\delta > 0$ is small enough so that for every $z \in H$ and every ${\widehat{z}} \in \D$ such that $|{\widehat{z}} - z| \le \delta$, the hyperbolic distance between $z$ and ${\widehat{z}}$ is at most $C_0/2$.

Let $\varepsilon > 0$ be such that every $K_1$-qc homeomorphism $\varphi : \D \to \D$ fixing $z = 0$ that is conformal outside a set of area $\varepsilon$, is at Euclidean distance at most $\delta$ form a map of the form $z \mapsto \lambda z$, with $|\lambda| = 1$ (part~2 of Lemma~\ref{l:approximation conformal}.)

Given $i \ge 0$ let $\varphi_i : \D \to \D$ be a M\"obius map such that $\varphi_i(H_i) = H$.
As $\| E \|_\D \le \varepsilon$ by hypothesis, it follows that $\chi \circ \varphi_i$ is conformal outside a set of Lebesgue measure at most $\varepsilon$.
Thus there is a M\"obius map $\psi_i : \D \to \D$ such that $\psi_i \circ \chi \circ \varphi_i (0) = 0$ and such that the Euclidean distance from $\psi_i \circ \chi \circ \varphi_i$ to the identity is at most $\delta$.

Thus, the hyperbolic distance between $\chi$ and $(\varphi_i \circ \psi_i)^{-1}$ on $\ov{H_i}$ is at most $C_0 /2$.
Moreover, it follows that there is a homeomorphism $\chi_i : \ov{H_i} \to \chi(\ov{H_i})$ that is $K_0$-qc on the interior of $H_i$ and whose hyperbolic distance to $(\varphi_i \circ \psi_i)^{-1}|_{\ov{H_i}}$ is at most $C_0/2$.
Thus the hyperbolic distance between $\chi_i$ and $\chi|_{\ov{H_i}}$ is at most $C_0$.

Is easy to see that the map $\hchi : \D \to \D$ that coincides with $\chi_i$ on $\ov{H_i}$, for every $i \ge 0$, is well defined and that it is a homeomorphism.
Moreover $\hchi$ is $K_0$-qc outside $\cup_{i \ge 0} \partial H_i$.
As this set has $\sigma$-finite length, it is qc removable, so $\hchi$ is a $K_0$-qc homeomorphism whose hyperbolic distance to $\chi$ is at most $C_0$.
\end{proof}

\subsection{Proof of Rigidity.}\label{ss:proof of Rigidity}
Let $R$ be a rational map of degree at least~$2$ and let $K > 1$ be given.
Let $\modul > 0$ be sufficiently large so that the constant $K(\modul)$ given by Lemma~\ref{l:smoothing3} is less than $K$ and put $K_0 \= KK(\modul)^{-1} > 1$ and $C_1 \= C(K_0, \modul)$.
Choose $C_0 > 0$ and let $\varepsilon > 0$ be given by Lemma~\ref{l:smoothing1} for $K_1 = K$ and for this choice of $K_0$ and $C_0$.

Let $(\hhV, \hV, V)$ be a nice triple for $R$, satisfying the hypothesis of Rigidity.
For each $c \in \CJ$ put
$$
A^c \= \hV^c \setminus \ov{V^c}
\text{ and }
A \= \cup_{c \in \CJ} A^c.
$$

Denote by $\cW$ the collection of connected components of $\hV \setminus K(V)$.
Note that for every $c \in \CJ$ we have $V^c \in \cW$.
For $W \in \cW$ the map $R^{m_W}$ in univalent on $W$ and the inverse of $R^{m_W}|_W$ extends univalently to $\hhV^{c(W)}$.
We denote this (extended) inverse by $\phi_W$ and put $\hW \= \phi_W(\hV^{c(W)})$.
Moreover, we denote by $\ell(W) \ge 0$ be the number of integers $m = 1, \ldots, m(U)$ for which $R^m(W) \subset \hV$.
Note that $\ell(W) = 0$ if and only if for some $c \in \CJ$ we have $W = V^c$.
Moreover, for every $W \in \cW$ with $\ell(W) \ge 1$ there is a (unique) $W' \in \cW$ such that $\ell(W') = \ell(W) - 1$ and such that $W$ intersects $\hW'$.
In this case $\hW \subset \hW'$.
\begin{lemm}\label{l:Rigidity induction step}
There is a constant $C_2 > 0$ only depending on
\begin{equation}\label{e:modulus and dilatation}
\min_{c \in \CJ} \modulus \left( \hhV^c \setminus \ov{\hV^c} \right)
\text{ and }
K,
\end{equation}
such that for every positive integer $\ell_0$ and every $c \in \CJ$ the following property holds.
Let $\chi : \hV^c \to \chi(\hV^c)$ be a qc map that is \clae{} on $(\hV^c \cap K(V)) \setminus (\cup_{W \in \cW, \ell(W) = \ell_0} \hW)$ and that is $K$-qc on $\hV^c \cap (\cup_{W \in \cW, \ell(W) = \ell_0} \hW)$.
Then there is a $K$-qc homeomorphism $\hchi : \hV^c \to \chi(\hV^c)$ whose hyperbolic distance to $\chi$ is at most $C_0 + C_1 + C_2$.
\end{lemm}
\begin{proof}
Denote by $\hcW$ (resp. $\hhcW$) the collection of connected components of $A \setminus R^{-1}(K(\hV))$ (resp. $A \setminus R^{-1}(K(\hhV))$).
For each element $W$ of $\hcW$ or $\hhcW$, the map $R^{m_W}$ is univalent on $W$ and the inverse of $R^{m_W}|_W$ extends univalently to $\hhV^{c(W)}$.
We will denote this (extended) inverse map by $\phi_W$.
Note that for $W \in \cW$ we have $\hW \in \hcW$ if and only if $\ell(W) = 1$.

For $c \in \CJ$ put $\hA^c \= \hhV^c \setminus \ov{\hV^c}$.
Note that for each $\hW \in \hcW$ contained in $\hV^c$, the set $\phi_{\hW}(\hA^{c(\hW)})$ is an annulus in $A^c$ that encloses $\hW$ and whose modulus is equal to that of $\hA^{c(\hW)}$.
It follows that there is a constant $C_2 > 0$, only depending on~\eqref{e:modulus and dilatation}, such that for every $K$-qc homeomorphism $\psi : A^c \to \psi(A^c)$ the hyperbolic diameter of $\psi(\hW)$ in $\psi(A^c)$ is at most $C_2$.
In particular note that if $\tchi : \hV^c \to \chi(\hV^c)$ is a homeomorphism that coincides with $\chi$ on $\hV^c \setminus \cup_{\hW \in \hcW} \hW$, then for every $\hW \in \hcW$ contained in $A^c$ we have $\tchi(\hW) = \chi(\hW)$ and therefore the hyperbolic distance between $\tchi$ and $\chi$ is at most~$C_2$.

\partn{1}
We will prove the assertion for $\ell_0 = 1$, with $C_0 + C_1 + C_2$ replaced by $C_0 + C_1$.
So by assumption $\chi$ is \clae{} on $A^c \setminus \cup_{\hW \in \hcW} \hW$ and $K$-qc on $A^c \cap (\cup_{\hW \in \hcW} \hW)$.

For each $\hhW \in \hhcW$ contained in $A^c$, the map $\chi \circ \phi_{\hhW} : \hhV^{c(\hhW)} \to \chi(\hhW)$ is a $K$-qc homeomorphism onto its image that is \clae{} on $\hhV^{c(\hhW)} \cap K(\hV)$.
Then Lemma~\ref{l:smoothing1} implies that there is a $K_0$-qc homeomorphism $\chi_{\hhW} : \hhW \to \chi(\hhW)$ whose hyperbolic distance to $\chi_1|_{\hhW}$ is at most $C_0$.
By the Gluing Lemma it follows that there is a qc homeomorphism $\chi_0 : \hV^c \to \chi(\hV^c)$ that coincides with $\chi$ on $\hV^c \setminus \cup_{\hhW \in \hhcW} \hhW$ and such that for every $\hhW \in \hhcW$ it coincides with $\chi_{\hhW}$ on $\hhW$.
It particular $\chi_0$ is $K_0$-qc on $A^c$ and the hyperbolic distance between $\chi_0$ and $\chi$ is at most $C_0$.
So Lemma~\ref{l:smoothing3} applied to $\chi = \chi_0$ implies that there is a $K_0K(\modul)$-qc homeomorphism $\hchi : \hV^c \to \chi(\hV^c)$ whose hyperbolic distance to $\chi_0$ is at most $C(K_0, \modul)$.
It follows that the hyperbolic distance between $\chi$ and $\hchi$ is at most $C_0 + C(K_0, \modul)$.
As by definition $K = K_0K(\modul)$ and $C_1 = C(K_0, \modul)$, this proves the assertion.

\partn{2}
Suppose that $\ell_0 > 1$ and assume by induction hypothesis that the assertion of the lemma holds for $\ell_0$ replaced by $\ell_0 - 1$.

Let $W \in \cW$ be satisfying $\ell(W) = 1$, so that $\hW \in \hcW$.
Then the map $\chi \circ \phi_W|_{\hV^{c(W)}}$ satisfies the hypothesis of the lemma with $\ell_0$ replaced by $\ell_0 - 1$.
By the induction hypothesis there is a $K$-qc map $\chi_W : \hW \to \chi(\hW)$ whose hyperbolic distance to $\chi|_W$ is finite.
By the Gluing Lemma it follows that there is a qc-map $\tchi$ that is $K$-qc on $A^c \cap (\cup_{\hW \in \hcW} \hW)$, \clae{} on $A^c \setminus \cup_{\hW \in \hcW} \hW$, and that coincides with $\chi$ on $\hV^c \setminus (\cup_{\hW \in \hcW} \hW)$.
It follows that the hyperbolic distance between $\chi$ and $\tchi$ is at most~$C_2$.

By part~$1$ it follows that there is a $K$-qc homeomorphism $\hchi : \hV^c \to \chi(\hV^c)$ whose hyperbolic distance to $\tchi$ is at most $C_0 + C_1$.
Thus the hyperbolic distance between $\chi$ and $\hchi$ is at most $C_0 + C_1 + C_2$.
This completes the proof of the lemma.
\end{proof}

\begin{lemm}\label{l:rigidity first step}
Let $c \in \CJ$ and let $\chi : \hV^c \to \chi(\hV^c)$ be a qc homeomorphism that is \clae{} on $\hV^c \cap K(V^c)$.
Then there is a sequence $(\chi_{\ell})_{\ell \ge 0}$ of qc homeomorphisms from $\hV^c$ to $\chi(\hV^c)$, such that $\chi_\ell$ is conformal on a set of full Lebesgue measure on the complement of $\cup_{W \in \cW, \ell(W) \le \ell} W$, and such that $\chi_\ell$ converges uniformly on compact sets to $\chi$ as $\ell \rightarrow \infty$.
\end{lemm}
\begin{proof}
Without loss of generality assume $\chi(\hV^c) = \D$ and $\chi(c) = 0$.
Let $\mu$ be the Beltrami differential of $\chi$ and for $\ell \ge 0$ let $\mu_\ell$ be the Beltrami differential that coincides with $\mu$ outside $\cup_{W \in \cW, \ell(W) > \ell} W$ and that is equal to~$0$ on this set.
By the Measurable Riemann Mapping Theorem it follows that there is a qc homeomorphism $\chi_\ell : \hV^c \to \D$ whose Beltrami differential is equal to $\mu_\ell$ and such that $\chi_\ell (c) = 0$.

By definition $\chi_\ell$ is conformal on a set of full Lebesgue measure outside $\cup_{W \in \cW, \ell(W) \le \ell} W$.
Moreover the map $\chi_\ell \circ \chi^{-1} : \D \to \D$ is \clae{} outside $\chi(\cup_{W \in \cW, \ell(W) > \ell} W)$ and its dilatation is bounded by that of $\chi$.
As the Lebesgue measure of this set converges to~$0$ as $\ell \rightarrow \infty$ (because the sets $\cup_{W \in \cW, \ell(W) = \ell} W$ are pairwise disjoint), it follows that every convergent subsequence of $(\chi_\ell \circ \chi^{-1})_{\ell \ge 0}$ converges to a conformal homeomorphism of $\D$ that fixes~$0$ (part~1 of Lemma~\ref{l:approximation conformal}.)
But the only conformal homeomorphisms of $\D$ that fix $0$ are of the form $z \mapsto \lambda z$, for some $\lambda \in \C$ satisfying $|\lambda| = 1$.
Since the collection of $K$-qc homeomorphisms of $\D$ that fix $0$ is compact, we can replace each $\chi_\ell$ by a map of the form $\lambda_\ell \chi_\ell$, with $|\lambda_\ell| = 1$, such that $\chi_\ell \circ \chi^{-1}$ converges uniformly to the identity as $\ell \to \infty$.
\end{proof}
\begin{proof}[End of the proof of Rigidity.]
Let $c \in \CJ$ and let $\chi : \hV^c \to \chi(\hV^c)$ be a qc homeomorphism that is \clae{} on $\hV^c \cap K(V)$.

Let $( \chi_\ell )_{\ell \ge 0}$ be given by Lemma~\ref{l:rigidity first step}.
Note that for each $\ell \ge 0$ the qc homeomorphism $\chi = \chi_\ell$ satisfies the hypothesis of Lemma~\ref{l:Rigidity induction step} with $\ell_0 = \ell + 1$.
Thus there is a $K$-qc homeomorphism $\hchi_\ell : \hV^c \to \chi(\hV^c)$ whose hyperbolic distance to $\chi_\ell$ is bounded by $C_0 + C_1 + C_2$.
As $\chi_\ell$ converges uniformly on compact sets to $\chi$ as $\ell \rightarrow \infty$, it follows that a subsequence of $(\hchi_\ell)_{\ell \ge 0}$ converges uniformly on compact sets to a $K$-qc homeomorphism $\hchi : \hV^c \to \chi(\hV^c)$ whose hyperbolic distance to $\chi$ is at most $C_0 + C_1 + C_2$.
This proves that the nice couple $(\hV, V)$ is $K$-rigid.
\end{proof}

\section{Thurston's algorithm.}\label{s:Thurston's algorithm}
In this section we a sufficient condition for Thurston's algorithm to converge in a very specific situation.
In~\S\ref{ss:Thurston's algorithm} we recall Thurston's algorithm and state the main result of this section.
The proof of this result is in~\S\S\ref{ss:dynamics of tR},~\ref{ss:limit geometry}.
\subsection{Thurston's algorithm.}\label{ss:Thurston's algorithm}
Let $\cB$ be a set of~$3$ distinct \textsf{\textit{base points}} of the sphere $S^2$ and let $\tR : S^2 \to S^2$ be a ramified covering.
Each  homeomorphism $h : S^2 \to \CC$ defines a unique complex structure on $S^2$, for which $h \circ \tR : S^2 \to \CC$ is holomorphic.
So there is a unique  homeomorphism $\hh : S^2 \to \CC$ that is holomorphic with respect to this complex structure and that coincides with $h$ on $\cB$.
It follows that $Q \= h \circ \tR \circ \hh^{- 1} : \CC \to \CC$ is a rational map of the same degree as $\tR$. 
In this way we associate to each  homeomorphism $h : S^2 \to \CC$ a  homeomorphism $\hh : S^2 \to \CC$ that coincides with $h$ on $\cB$ and a rational map $Q$ of the same degree as $\tR$, such that $Q \circ \hh \equiv h \circ \tR$.
By iterating this procedure, to each  homeomorphism $h_0 : S^2 \to \CC$ we associate a sequence $( h_k )_{k \ge 0}$ of  homeomorphisms from $S^2$ to $\CC$ and a sequence of rational maps $( Q_k )_{k \ge 0}$ of the same degree as $\tR$, such that for every $k \ge 0$ the homeomorphism $h_k$ coincides with $h_0$ on $\cB$ and such that we have $Q_k \circ h_{k + 1} = h_k \circ \tR$.

For our purposes we identify $S^2$ with $\CC$ and we will say that a homeomorphism of $\CC$ is \textsf{\textit{normalized}} if it fixed each element of $\cB$.
Furthermore, the map $\tR$ will be quasi-regular (a rational map post-composed with a qc homeomorphism) and $h_0$ will be the identity map of~$\CC$.
It follows that for every $k \ge 0$ the homeomorphism $h_k$ of $\CC$ is normalized and~qc.

Before the statement of the main result of this section, we consider the following definition.
\begin{defi}
Let $R$ and $Q$ be rational maps of degree at least~$2$ and let $V$ be a nice set for $R$.
We say that a qc homeomorphism $\chi : \CC \to \CC$ is a $V$~\textsf{\textit{pseudo-conjugacy}} between $R$ and $Q$, if $\chi \circ R \equiv Q \circ \chi$ on $\CC \setminus V$ and if $\chi$ is \clae{} on $K(V)$.
\end{defi}
\begin{theoalph}[Convergence of Thurston's algorithm]\label{t:Thurston's algorithm}
Let $R$ be a rational map of degree at least~$2$ and let $(\hhV, \hV, V)$ be a nice triple such that for some $K > 1$ the nice couple $(\hV, V)$ is $K$-rigid and such that for every $c \in \CJ$ the set $\hhV^c \cap K(\hV)$ has positive Lebesgue measure.
Let $\tR$ be a quasi-regular map of the same degree as $R$, that coincides with $R$ outside $V$, and such that the image by $\tR$ of every ramification point of $\tR$ in $V$ is contained in $K(\hhV)$.
Moreover, choose a set $\cB$ of~$3$ base points in $K(\hhV)$, let $h_0$ be the identity map of $\CC$ and let $( h_k )_{k \ge 0}$ and $( Q_k )_{k \ge 0}$ be given by Thurston's algorithm, as above.

Then $Q_k$ converges uniformly to a rational map $Q$ and $h_k$ converges uniformly to a continuous map $h : \CC \to \CC$, as $k \to \infty$, and we have,
$$
Q \circ h \equiv h \circ \tR
\text{ on } \CC.
$$
Furthermore, $h$ maps ramification points of $\tR$ to critical points of $Q$ preserving local degrees, and there is a $K$-qc $\hV$~pseudo-conjugacy between $R$ and $Q$ that coincides with $h$ on $K(\hV)$.
\end{theoalph}
\begin{rema}\label{r:Thurston's algorithm}
\

\begin{enumerate}
\item[1.]
In general the map $h$ is not injective.
For example, is easy to construct a quasi-regular perturbation $\tR$ as in Theorem~\ref{t:Thurston's algorithm}, having a saddle periodic point.
As rational maps cannot have saddle periodic points, in this case $\tR$ cannot be conjugated to $Q$ and the map $h$ cannot be injective.
\item[2.]
For $c \in \CJ$ the quasi-regular map $\tR$ may have several ramification points in $V^c$, but the sum of their multiplicities must be equal to the multiplicity of $c$ as a critical point of $R$. 
\item[3.]
The conclusions of Theorem~\ref{t:Thurston's algorithm} hold under the weaker assumption that there is an integer $N$ such that for every ramification point $r$ of $\tR$ we have  $\tR^N(r) \in K(\hhV)$.	
The proof of this fact is essentially the same as the proof of Theorem~\ref{t:Thurston's algorithm}.
With similar arguments is possible to allow $\tR$ to have recurrent ramification points, with a Topological Collet-Eckmann dynamics. 
\end{enumerate}
\end{rema}
\subsection{Dynamics of $\tR$.}\label{ss:dynamics of tR}
Let $R$, $(V, \hV, \hhV)$, $K$ and $\tR$ by as in the statement of Theorem~\ref{t:Thurston's algorithm}.
Recall that for a connected component $\hW$ of $\CC \setminus K(\hV)$, the integer $m_{\hW} \ge 0$ and $c(\hW) \in \CJ$ are such that $R^{m_{\hW}}$ maps $\hW$ univalently onto $\hV^{c(\hW)}$.

For $m \ge 0$ put
$$
\tK(m) \= \tR^{-m}(K(\hV)),
$$
which is a compact set that is forward invariant by $\tR$.
Let $\tW$ be a connected component of $\CC \setminus \tK(m)$.
Then $\hW \= \tR^m(\tW)$ is a connected component of $\CC \setminus K(\hV)$.
We put $m_{\tW} \= m + m_W$ and $c(\tW) \= c(W)$, so that $R^{m_{\tW} - m} \circ \tR^m$ maps $\tW$ univalently onto $\hV^{c(\tW)}$.
In the case when $\tW$ is also a connected component of $\CC \setminus K(\hV)$, is easy to see that this definition of $m_{\tW}$ and $c(\tW)$ is consistent with the previous definition.
Note that $\tW$ is not a connected component of $\tK(m_{\tW} + 1)$, but for every $m' \le m_{\tW}$ the set $\tW$ is a connected component of $\CC \setminus \tK(m')$.

As the image by $\tR$ of the ramification points of $\tR$ belong to $K(\hhV)$ and this set is forward invariant by $\tR$, it follows that $R^{m_{\tW} - m} \circ \tR^m$ maps $\tW$ univalently onto $\hV^{c(\tW)}$ and its inverse extends to $\hhV^{c(\tW)}$.
We denote by $\phi_{\tW}$ the inverse of this map and put
$$
\tW^- \=  \phi_{\tW} \left( \hhV^{c(\tW)} \right)
\text{ and }
A_{\tW}
\=
\phi_{\tW} \left( \hhV^{c(\tW)} \setminus \ov{\hV^{c(\tW)}} \right)
=
\tW^- \setminus \ov{\tW}.
$$
So $A_{\tW}$ is an annulus contained in $\tW^-$ which encloses $\tW$.
In particular $A_{\tW}$ is contained in the complement of $K(\hhV)$.

Note that for $k \ge m \ge 0$ the rational map $Q_m \circ \ldots \circ Q_{k - 1}$ maps $h_k(\tW^-)$ univalently onto $h_m(\tW^-)$.
Taking $m = 0$ we conclude that the rational map
$$
R^{m_{\tW} - m} \circ Q_0 \circ \ldots \circ Q_{m - 1}
$$
maps $h_m(\tW^-)$ univalently onto $\hhV^{c(\tW)}$.

\begin{prop}\label{p:Cauchy argument}
\

\begin{enumerate}
\item[1.]
For $m \ge 0$ and $k \ge m$ there is a $K$-qc homeomorphism $h_{k, m} : \CC \to \CC$ that coincides with $h_k \circ h_m^{-1}$ on $h_m(\tK(m))$ and that is \clae{} on this set.
\item[2.]
For each $m \ge 0$ there is a normalized $K$-qc homeomorphism $\chi_m : \CC \to \CC$ that is \clae{} on $h_m(\tK(m))$ and such that for every $k \ge m$ the homeomorphisms $\chi_k \circ h_k$ and $\chi_m \circ h_m$ coincide on $\tK(m)$.
\end{enumerate}
\end{prop}
\begin{proof}

\partn{1}
As $\tR$ is a quasi-regular map and $h_0$ is equal to the identity, it follows that $h_\ell$ is qc for all $\ell \ge 0$.
Since for every $\ell \ge 1$ we have $Q_0 \circ \cdots \circ Q_{\ell - 1} \circ h_\ell = \tR^\ell$, for $k \ge m$ we have
$$
Q_0 \circ \cdots \circ Q_{k - 1}  \circ (h_k \circ h_m^{-1})
=
\tR^k \circ h_m^{-1}
=
\tR^{k -m} \circ Q_0 \circ \cdots \circ Q_{m - 1}.
$$
As $\tR$ coincides with $R$ on $K(V)$ it follows that $\tR$ this set is forward invariant by $\tR$ and that $\tR$ is \clae{} on this set.
On the other hand we have,
$$
K(V) = Q_0 \circ \cdots \circ Q_{m - 1} \circ h_m \circ \tR^{-m}(K(V)),
$$
so $h_k \circ h_m^{-1}$ is \clae{} on $h_m(\tR^{-m}(K(V)))$.

Let $\tW$ be a connected component of $\CC \setminus \tK(m)$ and set $\hW \= \tR^m(\tW)$, which is a connected component of $\CC \setminus K(\hV)$.
Note that the map
$$
\varphi_{\tW} \= R^{m_W} \circ Q_0 \circ \cdots \circ Q_{m - 1} |_{h_m(\tW)}
$$
is a biholomorphism between $h_m(\tW)$ and $\hV^{c(W)}$ that maps points in $h_m(\tR^{-m}(K(V)))$ to points in $K(V)$.
Thus the map $h_k \circ h_m^{-1} \circ \varphi_{\tW}^{-1}$ is a qc homeomorphism between $\hV^{c(W)}$ and $h_k(\tW)$ that is \clae{} on $\hV^{c(W)} \cap K(V)$.
As the nice couple $(\hV, V)$ is $K$-rigid, there is a $K$-qc homeomorphism $\psi_{\tW} : h_m(\tW) \to h_k(\tW)$ at bounded hyperbolic distance from $h_k \circ h_m^{-1}|_{h_m(\tW)}$.

By the Gluing Lemma there is a $K$-qc homeomorphism $h_{k, m}$ that coincides with $h_k \circ h_m^{-1}$ on $h_m(\tK(m))$ and with $\psi_{\tW}$ on $h_m(\tW)$, for every connected component $\tW$ of $\CC \setminus \tK(m)$.
As $h_k \circ h_m^{-1}$ is \clae{} on $h_m(\tK(m)) \subset h_m(\tR^{-m}(K(V))$, it follows that $h_{k, m}$ is also \clae{} on this set.

\partn{2}
As for every $m \ge 0$ the homeomorphism $h_m$ is normalized, it follows that for every $k \ge m$ the homeomorphisms $h_k \circ h_m^{-1}$ and $h_{k, m}$ are normalized.
Since the collection of normalized $K$-qc homeomorphisms of $\CC$ is compact, it follows by a diagonal argument that there is a sequence $( \ell_j )_{j \ge 0}$ such that for every $m \ge 0$ the homeomorphisms $h_{\ell_j, m}$ converge to a normalized $K$-qc homeomorphism $\chi_m : \CC \to \CC$, as $j \rightarrow \infty$.
Note that for $\ell \ge k \ge m$ the maps $h_k \circ h_m^{-1}$ and $h_{k, \ell}^{-1}  \circ h_{\ell, m}$ coincide on $h_m(\tK(m))$.
Moreover it follows that $h_k \circ  h_m^{-1}$ and $\chi_k^{-1} \circ \chi_m$ coincide on $h_m(\tK(m))$.
\end{proof}
\subsection{Limit geometry.}\label{ss:limit geometry}
For $m \ge 0$ put $K(m) = \chi_m \circ h_m (\tK(m))$.
Given a \cc{} $W$ of $\CC \setminus K(m)$ let $\tW = (\chi_m \circ h_m)^{-1}(W)$ and define $m_W = m_{\tW}$, $c(W) = c(\tW)$ and $A_W = \chi_m \circ h_m (A_{\tW})$.

For $c \in \CJ$ put $A^c \= \hhV^c \setminus \ov{\hV^c}$.
\begin{lemm}\label{l:image annuli}
For every $m \ge 0$ and every \cc{} $W$ of $\CC \setminus K(m)$ there is a $K$-qc homeomorphism $\psi_W : A^{c(W)} \to A_{W}$ that maps $A^{c(W)} \cap K(\hV)$ onto $A_W \cap K(m)$.
In particular there is a constant $\modul > 0$ independent of $m$ and $W$, such that $\modulus (A_W) \ge \modul$.
\end{lemm}
\begin{proof}
Put $\tW = (\chi_m \circ h_m)^{-1}(W)$ and note that the map $\varphi_{W} = R^{m_{W} - m} \circ Q_0 \circ \ldots \circ Q_{m - 1}$ is a biholomorphism between $h_m(\tW^-)$ and $\hhV^{c(W)}$ that maps points in $h_m(\tK(m))$ to points in $K(\hV)$.
As $\chi_m$ is $K$-qc and $h_m (A_{\tW}) = \varphi_{W}^{-1} (A^{c(W)})$, the map $\psi_W = \chi_m \circ \varphi_{W}^{-1}|_{A^{c(W)}}$ satisfies the desired properties.
\end{proof}
\begin{prop}\label{p:limit geometry}
For $m \ge 0$ put $K(m) = \chi_m \circ h_m (\tK(m))$.
\begin{enumerate}
\item[1.]
For every $\varepsilon > 0$ there is $M \ge 0$ such that for every $m \ge M$, the diameter of each connected component of $\CC \setminus K(m)$ is at most $\varepsilon$.
\item[2.]
The set $K(m)$ is increasing with~$m$ and $\cup_{m \ge 0} K(m)$ has full Lebesgue measure in $\CC$.
\end{enumerate}
\end{prop}
\begin{proof}
\

\partn{1}
Let $m \ge 0$ and let $W$ be a connected component of $\CC \setminus K(m)$.
Thus $\tW = (\chi_m \circ h_m )^{-1}(W)$ is a connected component of $\CC \setminus \tK(m)$.

Let $0 \le m_1 < \ldots < m_\ell < m_{\ell + 1} = m_{\tW}$ be all the integers $j = 0, \ldots, m_{\tW}$ such that $\tR^j(\tW)$ is contained in $\hV$.
Note that $m_\ell < m$.
For $i = 1, \ldots \ell + 1$ let $c_i \in \CJ$ be such that $\tR^{m_i}(\tW) \subset \hV^{c_i}$, so $c_{\ell + 1} = c(\tW)$.
Moreover denote by $\tW_i$ the connected component of $\tR^{- m_i}(\hV^{c_i})$ that contains $\tW$ and set $\tA_i = A_{\tW_i}$.
Note that $\tW_{\ell + 1} = \tW$, $\tW_1$ is a connected component of $\CC \setminus K(\hV)$ and $m_1 = m_{\tW_1}$.
Moreover, for $i = 1, \ldots, \ell$ we have $\tW_{i + 1}^- \subset \tW_i$, so the annuli $\tA_1, \ldots, \tA_{\ell + 1}$ are pairwise disjoint and enclose $\tW_{\ell + 1} = \tW$.

When $m$ is big there are three cases.

\partn{Case 1}
$m_1$ is big.
As $\tU = \tW_1$ is a connected component of $\CC \setminus K(\hV)$ and $m_1 = m_{\tU}$, it follows that $R^{m_{\tU}}$ maps $\tU$ biholomorphically onto $\hV^{c(\tU)}$.
Since the inverse of this map extends to $\hhV^{c(U)}$, a Fatou type argument shows that the diameter of $\tU$ is small when $m_{\tU}$ is big.
Since $\tW \subset \tU$, we have $W \subset \chi_0 \circ h_0 (\tU)$, so $W$ is small when $m_1$ is big.

\partn{Case 2}
There is $i \in \{ 1, \ldots, \ell + 1 \}$ such that $m_{i + 1} - m_i$ is big.
Put $\tU = \tR^{m_i}(\tW_{i + 1}) \subset \hV^{c_i}$ and $U = \chi_m \circ h_m (\tU)$.
In this case $\tR(\tU)$ is a connected component of $\CC \setminus K(\hV)$ and $m_{R(\tU)} = m_{i + 1} - m_i - 1$, so as in the previous case, the diameter of $\tU$ is small as $m_{i + 1} - m_i$ is big.
Thus the modulus of the annulus $\tA = \hhV^{c_i} \setminus \ov{\tU}$ is big as $m_{i + 1} - m_i$ is big.

As $Q_1 \circ \ldots \circ Q_m$ maps $h_{m + 1}(\tW_i^-)$ biholomorphically to $h_1(\hhV^{c_i})$ and $\chi_{m + 1}$ is a $K$-qc map, we have
$$
\modulus (\chi_{m + 1} \circ h_{m + 1}(\tW_i^- \setminus \ov{\tW_{i + 1}}))
\ge
K^{-1} \modulus (h_1(\tA)).
$$
Since this annulus is contained in the complement of $K(0)$, it follows that $W \subset \chi_{m + 1} \circ h_{m + 1} (\tW_{i + 1})$ is small as $m_{i + 1} - m_i$ is big.

\partn{Case 3}
$\ell$ is big.
Set $U = \chi_m \circ h_m (\tW_1^-)$.
Note that the annuli $\chi_m \circ h_m(A_1), \ldots, \chi_m \circ h_m(A_{\ell + 1})$ are pairwise disjoint and enclose $W$.
So, if $\modul > 0$ is the constant given by Lemma~\ref{l:image annuli}, Gr\"otzsch inequality implies that $\modulus(U \setminus \ov{W}) \ge \ell \modul$.
As $U$ is contained in the complement of $K(0)$, it follows that the diameter of $W$ is small as $\ell$ is big.

\partn{2}
By part~$2$ of Proposition~\ref{p:Cauchy argument}, for every $k \ge m$ the maps $\chi_k \circ h_k$ and $\chi_m \circ h_m$ coincide on $\tK(m)$.
As the set $\tK(m)$ is increasing with~$m$, it follows that~$K(m) = \chi_m \circ h_m (\tK(m))$ is also increasing with~$m$.

To prove that $\cup_{m \ge 0} K(m)$ has full Lebesgue measure on $\CC$, consider a point $p$ in the complement of this set.
With a reasoning similar to that of part~$1$ we can find an infinite sequence of pairwise disjoint annuli $A_1, A_2, \ldots$ that enclose $p$ and that are contained in the complement of $K(0)$.
By Lemma~\ref{l:image annuli} each $j \ge 1$ there are $m_j \ge 0$, $c_j \in \CJ$ and a $K$-qc homeomorphism $\psi_j : A^{c_j} \to A_j$ that maps $A^{c_j} \cap K(\hV)$ to $A_j \cap K(m_j)$.
By hypothesis, for each $c \in \CJ$ the set $A^c \cap K(\hV) = \hhV^c \cap K(\hV)$ as positive Lebesgue measure.
By the area distortion property of $K$-qc maps (see e.g.~[A]), it follows that $p$ is not a density point of the complement of $\cup_{m \ge 0} K(m)$.
As $p$ was an arbitrary point in the complement of $\cup_{m \ge 0} K(m)$, it follows that this set has full Lebesgue measure.
\end{proof}

\begin{proof}[Proof of Theorem~\ref{t:Thurston's algorithm}]
\

\partn{1}
By part~$1$ of Proposition~\ref{p:Cauchy argument}, for $k \ge m$ the homeomorphisms $\chi_k \circ h_k$ and $\chi_m \circ h_m$ coincide on $\tK(m)$, so for each connected component $\tW$ of $\CC \setminus \tK(m)$ we have $\chi_k \circ h_k (\tW) = \chi_m \circ h_m (\tW)$.
As the diameter of $\chi_m \circ h_m(\tW)$ is small as $m$ is big (part 1 of Proposition~\ref{p:limit geometry}) it follows that $\chi_m \circ h_m$ converges uniformly to a continuous function $h : \CC \to \CC$.
Similarly the we can show that the map $\chi_m \circ Q_m \circ \chi_{m - 1}^{-1}$ converges uniformly to some map $Q : \CC \to \CC$.

Recall that $\chi_m$ is a normalized $K$-qc homeomorphism that is conformal on a set of full Lebesgue measure in $h_m(\tK(m))$ (part~2 of Proposition~\ref{p:Cauchy argument}.)
So $\chi_m^{- 1}$ is \clae{} on $K(m) = \chi_m \circ h_m (\tK(m))$.
Since $\cup_{k \ge 0} K(m)$ has full Lebesgue measure in $\CC$ (part~2 of Proposition~\ref{p:limit geometry}), it follows that every convergent subsequence of $( \chi_m^{-1} )_{m \ge 0}$ converges uniformly to a normalized conformal map as $m \rightarrow \infty$ (Lemma~\ref{l:approximation conformal}.)
But the only normalized conformal map is the identity, so every convergent subsequence of $( \chi_m^{-1} )_{m \ge 0}$ converges to the identity.
Since the collection of normalized $K$-qc homeomorphisms is compact, it follows $\chi_m^{-1}$, and hence $\chi_m$, converges to the identity as $m \to \infty$.

Hence we conclude that $h_m$ and $Q_m$ converge uniformly to $h$ and $Q$, respectively, as $m \rightarrow \infty$.
As the rational maps~$Q_m$ have the same degree as $R$, it follows that $Q$ is a rational map of the same degree as~$R$.

Let us prove now that $h$ maps ramification points of $\tR$ to critical points of $Q$, preserving local degrees.
First we will prove that $h$ does not decrease local degrees of points in $\tK(1) = \tR^{-1}(K(\hV))$.
Note that $h$ is injective on $\tK(1)$, as it coincides with the homeomorphism $\chi_1 \circ h_1$ there.

Given $\tp \in \tK(1)$ choose $\tq$ close enough to $h(\tp)$, so that the number of preimages of $\tq$ by $h$ near $\tp$ is equal to the local degree of $\tR$ at $\tp$.
As $h(\tp) \in K(\hV)$ and the set $K(\hV)$ is connected (part~2 of Lemma~\ref{l:nice sets forward}), we may choose $\tq$ in $K(\hV)$.
Thus, all preimages of $\tq$ by $\tR$ belong to $\tK(1) = \tR^{-1}(K(\hV))$.
As $h$ is injective on $\tK(1)$ and $K(\hV) \subset \tK(1)$, it follows that the local degree of $Q$ at $h(\tp)$ is at least as big as the local degree of $\tR$ at $\tp$.

Note that by hypothesis all ramification points of $\tR$ belong to $\tK(1)$.
Since $h$ is injective on $\tK(1)$, we conclude that the ramification points of $\tR$ are mapped to distinct critical points of $Q$, not decreasing multiplicities.
Since $\tR$ and $Q$ have the same degree, they have the same number of critical points counted with multiplicity, so we conclude that $h$ preserves the local degree of each ramification point.

\partn{2}
By part~2 of Proposition~\ref{p:Cauchy argument} the map $\chi_0$ is $K$-qc and it is \clae{} on $K(\hV)$.
Moreover, as $h_0$ is the identity, $\chi_0$ coincides with $h$ on $K(\hV)$.
It follows that $\chi_0$ conjugates $\tR$ to $Q$ on $K(\hV)$.
As $\tR$ and $R$ coincide on $K(\hV)$, it follows that $\chi_0$ conjugates $R$ and $Q$ on $K(\hV)$.

Let $\tW$ be a connected component of $\CC \setminus K(\hV)$ and set $W = \chi_0(\tW)$.
Then $R^{m_{\tW}}$ maps $\tW$ univalently to $\hV^{c(\tW)}$ and $Q^{m_{\tW}}$ maps $W$ univalently to $\chi_0(\hV^{c(\tW)})$.
Denote by $\chi_{\tW} : \tW \to W$ the $K$-qc homeomorphism $(Q^{m_{\tW}}|_W)^{-1} \circ \chi_0 \circ R^{m_{\tW}}|_{\tW}$.
It extends continuously to the closure of $\tW$ and it coincides with $R^{m_{\tW}}$ on the boundary of $\tW$.
Thus the Gluing Lemma implies that there is a $K$-qc homeomorphism $\chi$ that coincides with $h$ on $K(\hV)$ and with $\chi_{\tW}$ on $\tW$, for every connected component $\tW$ of the complement of $K(\hV)$.
By construction $\chi$ conjugates $R$ to $Q$ outside $\hV$.
Moreover $\chi$ coincides with $\chi_0$ and $h$ on $K(\hV)$ and it is \clae{} on $K(\hV)$.
\end{proof}

\section{\BCn{} and the \UPC.}\label{s:selfsimilar}
In this section we construct nice sets and nice nests for \BCg{} rational maps.
We begin by giving an equivalent formulation of the \BCn{} condition in \S\ref{ss:upc}.
After some preliminary lemmas in \S\ref{ss:lemmas}, we construct nice sets and nests in~\S\ref{ss:construction nice sets and nests}.

Throughout all this section we fix a rational map $R$ of degree at least~$2$.
\subsection{\UPC}\label{ss:upc}
Given $\delta' > \delta > 0$ we will say that a rational map $R$ of degree at least~$2$ satisfies the \textsf{\textit{\UP{$\delta$}{$\delta'$}}}, if for every $z \in \CC$ and every integer $n \ge 1$ such that,
\begin{enumerate}
\item[1.]
for $j = 1, \ldots n - 1$, we have $R^j(z) \not \in \tB(\CJ, \delta)$ and;
\item[2.]
for some $c \in \CJ$ we have $R^n(z) \in \tB(c, \delta')$;
\end{enumerate}
the pull-back of $\tB(c,\delta')$ to $z$ by $R^n$ is univalent.

Taking $\zeta = R(z)$ and $m = n -1$, we have the following equi\-va\-lent formulation of the \UP{$\delta$}{$\delta'$}:
If for $j = 0, \ldots, m - 1$ we have $R^j(\zeta) \not \in B(\CVJ, \delta)$ and for some $c \in \CJ$ we have $R^m(\zeta) \in \tB(c, \delta')$, then the pull-back of $\tB(c, \delta')$ to $\zeta$ by $R^m$ is univalent and disjoint from $\CV(R)$.

Note that if $\delta' > \delta > 0$ are such that $R$ satisfies the \UP{$\delta$}{$\delta'$}, then for every $\delta_0, \delta_0' \in (\delta, \delta')$ such that $\delta_0' > \delta_0$, the rational map $R$ satisfies the \UP{$\delta_0$}{$\delta_0'$}.
On the other hand, note that if for some $\delta' > \delta > 0$ the rational map $R$ is \BC{$\delta$}{$\delta'$}, then for every $\delta'' \in (\delta, \delta')$ the rational map $R$ is \BC{$\delta$}{$\delta''$}.
\begin{prop}\label{p:BC-UP}
The following assertions hold for a rational map $R$ of degree at least~$2$.
\begin{enumerate}
\item[1.]
For each $\delta' > \delta > 0$ the $(\delta, \delta')$\nobreakdash-\BCn{} condition implies the \UP{$\delta$}{$\delta'$}.
\item[2.]
There is a constant $\kappa \in (0, 1)$ only depending on~$R$, such that for every $\delta > 0$ and $\delta' > \kappa^{-1} \delta$ small, the \UP{$\delta$}{$\delta'$} implies the $(\delta, \kappa \delta')$\nobreakdash-\BCn{} condition.
\item[3.]
If $\kappa \in (0, 1)$ is as in part~$2$, then for every $\delta > 0$, $\delta' > \kappa^{-1}\delta$ and $\delta_0, \delta_0' \in (\delta, \kappa \delta')$ satisfying $\delta_0' > \delta_0$, the $(\delta, \delta')$\nobreakdash-\BCn{} condition implies the $(\delta_0,\delta_0')$\nobreakdash-\BCn.
\end{enumerate}
\end{prop}
\begin{proof}
\

\partn{1}
Let $z \in \CC$ and $n \ge 0$ be as in the \UPC{} and consider successive pull-backs $W_0, \ldots, W_n = \tB(c, \delta')$ by $R$ along the orbit of $z$, so that $z \in W_0$.
If for some $j = 0, \ldots, n - 1$ we have $W_j \cap \Crit \neq \emptyset$, then by the $(\delta, \delta')$\nobreakdash-\BCn{} condition, $\diam(W_{j + 1}) < \delta$.
Since by hypothesis $R^{j + 1}(z) \in W_{j + 1}$ does not belong to $B(\CV, \delta)$ we have that $W_{j + 1} \cap \CV = \emptyset$, which contradicts our assumption.
So the for every $j = 0, \ldots, n - 1$ the set $W_j$ is disjoint from $\Crit$ and therefore $R^n$ is univalent on $W_0$.

\partn{2}
Let $\kappa_0, \kappa_1 \in (0,1)$ to be determined below and assume that $\delta \le \kappa_0\kappa_1 \delta'$.
Given $c \in \CJ$ consider successive pull-backs $W_0, \ldots, W_m = \tB(c, \kappa_1 \kappa_0 \delta')$ by $R$, such that for $j = 0, \ldots, m - 1$ the set $W_j$ is disjoint from $B(\CV, \delta) = \emptyset$ and such that for some critical value $v \in \CVJ$ we have $W_0 \cap B(v, \delta) \neq \emptyset$.
We will prove that for appropriated choices of the constants $\kappa_0$ and $\kappa_1$, only depending on $R$, we have $\diam(W_0) < \delta \le \kappa_0 \kappa_1 \delta'$.
This implies by induction that $R$ is \BC{$\delta$}{$\kappa_0 \kappa_1 \delta'$} and proves the assertion of the proposition with constant $\kappa = \kappa_0 \kappa_1$.

For $j = 0, \ldots, m$ let $\hW_j$ (resp. $W_j'$) be the connected component of 
$$
R^{j - m}(\tB(c, \kappa_0 \delta'))
\text{ (resp. $R^{j - m}(\tB(c, \delta'))$)}
$$
that contains $W_j$, so that $W_j \subset \hW_j \subset W_j'$.
By the \UP{$\delta$}{$\delta'$} $R^m$ is univalent on $W_0'$.

\partn{Case 1}
$\hW_0 \subset B(v, \delta)$.
Since $R^m$ is univalent on $W_0'$, it follows by \KDT{} that the distortion of $R^m$ on $\hW_0$ can be bounded in terms of $\kappa_0$ only.
Thus given $\kappa_0 > 0$ we may choose $\kappa_1 \in (0, \tfrac{1}{2})$ small enough, such that
$$
\diam(W_0) < \tfrac{1}{2} \diam(\hW_0) \le \delta.
$$

\partn{Case 2} $\hW_0 \not \subset B(v, \delta)$.
By the equivalent formulation of the \UPC{} it follows that $W_0' \cap CV = \emptyset$. 
As $R^m$ is univalent on $W_0'$, the modulus of the annulus $W_0' \setminus \ov{\hW_0}$ is bounded from below by a function depending on $\kappa_0 \in (0, 1)$ only and it is big as $\kappa_0 > 0$ is small
(note that by Lemma~\ref{l:backward size} we have $\diam(\CC \setminus W_0') > \tfrac{1}{2} \diam(\CC)$.)
Since $\hW_0 \cap B(v, \delta) \neq \emptyset$ and $v \not \in W_0'$, we can choose $\kappa_0$ small enough, only depending on $R$, such that $\diam(\hW_0) < \delta$.

\partn{3} A direct consequence of parts~$1$ and~$2$.
\end{proof}
\subsection{Lemmas.}\label{ss:lemmas}
\begin{lemm}\label{l:prenice}
Suppose that for some $\tdelta > 0$ and $\tdelta' \ge 2 \tdelta$ the rational map $R$ is \BC{$\tdelta$}{$\tdelta'$}.
For $c \in \CJ$ let $\tV^c$ be the connected component of $\CC \setminus K(\tB(\CJ, \tdelta))$ that contains~$c$.
Then
$$
\tB(c, \tdelta) \subset \tV^c \subset \tB(c, 2\tdelta).
$$
\end{lemm}
\begin{proof}
Given $c \in \CJ$ and an integer $n \ge 0$, let $\tV_n^c$ be the connected component of $\cup_{j = 0, \ldots, n - 1} R^{-j}(\tB(\CJ, \tdelta))$ that contains $c$.
Note that $\tV_0^c = \tB(c, \tdelta)$, that $\tV_n^c$ is increasing with $n$, and that $\tV^c = \cup_{n \ge 0} \tV_n^c$.
So it is enough to show that for every $c \in \CJ$ and every integer $n \ge 0$ we have $\tV_n^{c} \subset \tB(c, 2\tdelta)$.

We will proceed by induction in~$n$.
For $n = 0$ the assertion is trivial.
Suppose that the assertion holds for $n \ge 0$ and fix $c \in \CJ$.
For every point $z \in \tV_{n + 1}^{c}$ there is an integer $m \in \{ 0, \ldots, n + 1 \}$ and $c_0 \in \CJ$ such that $R^m(z) \in \tB(c_0, \tdelta)$; let $m(z)$ be the least of such integers.
Let $X$ be a connected component of $\tV_{n + 1}^{c} \setminus \tB(c, \tdelta)$ and let $z \in X$ for which $m(z)$ is minimal among points in $X$.
Let $c_0 \in \CJ$ be such that $R^{m(z)}(z) \in \tB(c_0, \tdelta)$.
Considering that $m(z) > 0$, we have by induction hypothesis 
$$
R^{m(z)}(X) \subset \tV_n^{c_0} \subset \tB(c_0, 2\tdelta).
$$
As $R$ is \BC{$\tdelta$}{$\tdelta'$} it follows that $\diam(R(X)) < \tdelta$.
Therefore $R(X) \subset B(R(c), 2\tdelta)$ and $X \subset \tB(c, 2\tdelta)$.
This shows that $\tV_{n + 1}^c \subset \tB(c, 2 \tdelta)$ and completes the proof of the lemma.
\end{proof}
\begin{lemm}\label{l:diameter estimate}
There is a constant $C_0 > 0$ only depending on $R$ such that for each $\varepsilon \in (0, \tfrac{1}{2}]$ the following property holds.
Let $\delta' > \delta > 0$ be small and suppose that $R$ is \BC{$\delta$}{$\delta'$}.
Moreover, let $z \in \CC$ and $n \ge 1$ be such that for $j = 0, \ldots, n - 1$ we have $R^j(z) \not\in \tB(\CJ, \delta)$ and such that for some $c \in \CJ$ we have $R^n(z) \in \tB(c, \varepsilon\delta')$.
Then the pull-back $W$ of $B(c, \varepsilon \delta')$ by $R^n$ satisfies,
$$
\diam(W) \le C_0 \varepsilon^{\frac{1}{\mu_c}} \max \left\{ \delta, \dist(W, \CV) \right\}.
$$
\end{lemm}
\begin{proof}
Let $r_K > 0$ be the constant as in~\S\ref{ss:distortion}.
We assume that $\delta' > 0$ is small enough so that $\diam(\tB(c, \delta')) \le r_K$.

The \UPC{} implies that the pull-back $\hW$ of $\tB(c, \delta')$ by $R^n$ containing $z$ is univalent.
There are two cases.

\partn{Case 1} $\dist(W, \CV) < \delta$. 
As $R$ is \BC{$\delta$}{$\delta'$}, we have $\diam(\hW) < \delta$.
So there is a constant $C_0' > 0$ only depending on $R$ such that $\diam(W) \le C_0'\varepsilon^{\frac{1}{\mu_c}} \delta$.

\partn{Case 2} $\dist(z, \CV) \ge \delta$.
By the \UPC{} $\hW \cap \CV = \emptyset$.
As $\diam(\CC \setminus W) > \tfrac{1}{2}\diam(\CC)$ by Lemma~\ref{l:backward size}, there is a $C_0'' > 0$ only depending on $R$ such that $\diam(W) \le C_0'' \varepsilon^{\frac{1}{\mu_c}} \dist(W, \CV)$.
\end{proof}
\begin{lemm}\label{l:pre construction}
Let $R$ be a rational map of degree at least~$2$ and let $\teta > 1$ be given.
Let $\tV_0$ be a neighborhood of $\CJ$ disjoint from the forward orbits of critical points of $R$ not in $J(R)$ and such that every connected component of $\tV_0$ contains precisely one element of $\CJ$.
Suppose moreover that there is $\tdelta > 0$ such that for each $c \in \CJ$ the connected component $\tV^c$ of $\CC \setminus K(\tV_0)$ satisfies
$$
\tB(c, \tdelta) \subset \tV^c \subset \tB(c, \teta \tdelta).
$$
Let $C_0 > 0$ be the constant given by Lemma~\ref{l:diameter estimate} and put
$$
\eta = 1 + C_0(\teta\tdelta/\delta')^{\frac{1}{\mu_{\max}}}.
$$
Suppose that there are $\delta < \tdelta$ and $\delta' \ge 2\teta\tdelta$ small such that is $R$ \BC{$\delta$}{$\delta'$}.
Then for every $\hdelta \ge \delta$ and every $v \in \CVJ$ there is a simply-connected neighborhood $U^v$ of $v$ in $\CC$ such that $\partial U^v \subset K(\tV_0)$ and such that
$$
B(v, \hdelta) \subset U^v \subset B(v, \eta\hdelta).
$$
If moreover $\hdelta < \delta'$ and $R$ is \BC{$\hdelta$}{$\delta'$}, then we also have $U^v \subset B(v, 2\hdelta)$.
\end{lemm}
\begin{proof}
Let $W_0$ be a connected component of $\CC \setminus K(\tV_0)$ and let $m \ge 0$ and $c \in \CJ$ be such that $R^m(W_0) = \tV^c$.
The \UPC{} implies that the inverse of $R^m|_{W_0}$ extends in a univalent way to $B(c, \delta')$.
So Lemma~\ref{l:diameter estimate} implies that, if $W_0$ intersects $B(\CV, \hdelta)$, then
$$
\diam(W_0)
\le
C_0(\teta\tdelta/\delta')^{\frac{1}{\mu_{\max}}} \hdelta
=
(\eta - 1)\hdelta.
$$
 
Given $v \in \CVJ$ put,
\begin{multline*}
\widehat{U}^v
=
B(v, \hdelta) \cup \left\{ \text{ connected components of } \CC \setminus K(\tV_0)
\right.
\\
\left. \text{intersecting } B(v, \hdelta) \right\},
\end{multline*}
so that $B(v, \hdelta) \subset \widehat{U}^v \subset B(v, \eta\hdelta)$ and $\partial \widehat{U}^c \subset K(\tV_0)$.
Moreover let 
\begin{multline*}
U^v
=
\widehat{U}^v \cup \left\{ \text{ connected components $W_0$ of } \CC \setminus \widehat{U}^v \right.
\\
\left. \text{ such that } \diam(W_0) < \diam(\widehat{U}^v) \right\}.
\end{multline*}
As $\diam(\widehat{U}^v) \ll \diam(\CC)$, there is only one connected component of $\CC \setminus \widehat{U}^v$ whose diameter is not smaller than $\widehat{U}^v$.\
Thus $U^v$ is simply-connected.
Moreover $B(v, \hdelta) \subset U^v \subset B(v, \eta\hdelta)$ and $\partial U^v \subset \partial \widehat{U}^v \subset K(\tV_0)$.

If $\hdelta < \delta'$ and $R$ is \BC{$\hdelta$}{$\delta'$}, then every connected component $W_0$ of $\CC \setminus K(\tV_0)$ intersecting $B(\CVJ, \hdelta)$ satisfies $\diam(W_0) < \hdelta$.
So in this case we have $\hU^v \subset B(v, 2\hdelta)$ and hence $U^v \subset B(v, 2\hdelta)$.
\end{proof}
\subsection{Construction of nice sets and nests.}\label{ss:construction nice sets and nests}
We will say that a nice set $V = \cup_{c \in \CJ} V^c$ is \textsf{\textit{symmetric}} if for each $v \in \CVJ$ there is a simply-connected neighborhood $U^v$ of $v$ in $\CC$ such that for every $c \in \CJ$ the set $V^c$ is equal to the connected component of $R^{-1}(U^{R(c)})$ containing $c$.
\begin{prop}\label{p:nice set}
Let $R$ be a rational map of degree at least~$2$, let $\kappa \in (0, 1)$ be the constant given by Proposition~\ref{p:BC-UP} and let $C_0$ be the constant given by Lemma~\ref{l:diameter estimate}.
Suppose that for some $\delta > 0$ and $\delta' > 8 \kappa^{-1} \delta$ the rational map $R$ is \BC{$\delta$}{$\delta'$} and put
$$
\eta \= 1 + \min \left\{ 1, C_0(4\delta/\delta')^{\frac{1}{\mu_{\max}}} \right\}.
$$
Then there is a symmetric nice set $V = \cup_{c \in \CJ} V^c$ such that for every $c \in \CJ$ we have
$$
\tB(c, \delta) \subset V^c \subset \tB(c, \eta \delta).
$$
\end{prop}
\begin{proof}
Choose $\tdelta$ in $(2\delta, \min \{ 4\delta, \tfrac{1}{4}\kappa \delta' \})$.
Part~$3$ of Proposition~\ref{p:BC-UP} implies that $R$ is \BC{$\tdelta$}{$2\tdelta$}.
Lemma~\ref{l:prenice} implies that for each $c \in \CJ$ the connected component $\tV^c$ of $\CC \setminus K(B(\CJ, \tdelta))$ satisfies,
$$
\tB(c, \tdelta) \subset \tV^c \subset \tB(c, 2\tdelta).
$$
Lemma~\ref{l:pre construction} applied to $\teta = 2$ and $\hdelta = \delta$, implies that for each $v \in \CVJ$ there is a simply-connected neighborhood $U^v$ of $v$ in $\CC$ such that $\partial U^v \subset K(B(\CJ, \tdelta))$ and such that $B(v, \delta) \subset U^v \subset B(v, \eta\delta)$.
For each $c \in \CJ$ we denote by $V^c$ the connected component of $R^{-1}(V^{R(c)})$, so that $\tB(c, \delta) \subset V^c \subset \tB(c, \eta\delta)$.
As $\eta\delta \le 2\delta < \tdelta$, Lemma~\ref{l:nice sets with closures} implies that $V = \cup_{c \in \CJ} V^c$ is a nice set.
\end{proof}
\begin{prop}\label{p:nice nest}
Let $R$ be a rational map of degree at least~$2$ and let $C_0 > 0$ be the constant given by Lemma~\ref{l:diameter estimate}.
Moreover, let $\ell$ be a positive integer, $\tau \in (0, 1)$, $\eta \in (1, \min\{2, \tau^{-1}\})$ be given, and let $\dconst_0 > 1$ be such that
$$
\dconst_0 > 4\tau^{-(\ell + 1)} 
\text{ and }
C_0 \left( \tfrac{1}{2}\tau^{\ell + 1}\dconst_0 \right)^{-\frac{1}{\mu_{\max}}} \le \eta - 1.
$$
Then the following properties hold.
\begin{enumerate}
\item[1.]
Let $\kappa \in (0, 1)$ be the constant given by Proposition~\ref{p:BC-UP}.
Suppose that for some $\delta > 0$ the rational map $R$ is \BC{$\delta$}{$\kappa^{-1}\dconst_0\delta$} and put $\delta_0 = \tau^{-\ell}\delta$.
Then for each $j \in \{0, \ldots, \ell \}$ there is a symmetric nice set $V_j = \cup_{c \in \CJ} V_j^c$ such that for each $c \in \CJ$ we have,
$$
\tB(c, \tau^j\delta_0) \subset V_j^c \subset \tB(c, \eta\tau^j\delta_0),
$$
and such that the sequence $(V_0, \ldots, V_\ell)$ is a nice nest.
\item[2.]
Suppose that there is $\delta_0 > 0$ such that for every $\delta \in (0, \delta_0]$ the rational map $R$ is \BC{$\delta$}{$\dconst_0\delta$}.
Then for every positive integer $j \ge 1$ there is a symmetric nice set $V_j = \cup_{c \in \CJ} V_j^c$ such that for each $c \in \CJ$ we have,
$$
\tB(c, \tau^j\delta_0) \subset V_j^c \subset \tB(c, \eta\tau^j\delta_0),
$$
and such that for every positive integer $j_0 \ge 1$ the sequence $(V_{j_0}, \ldots, V_{j_0 + \ell})$ is a nice nest.
\end{enumerate}
\end{prop}
\begin{proof}
\

\partn{1}
Part~$3$ of Proposition~\ref{p:BC-UP} implies that $R$ is \BC{$\tau^{-1}\delta_0$}{$2\tau^{-1}\delta_0$}.
Lemma~\ref{l:prenice} applied to $\tau^{-1}\delta_0$ instead of $\tdelta$, implies that for each $c \in \CJ$ the connected component $\tV^c$ of $\CC \setminus K(B(\CJ, \tau^{-1}\delta_0))$ satisfies,
$$
\tB(c, \tau^{-1}\delta_0) \subset \tV^c \subset \tB(c, 2\tau^{-1}\delta_0).
$$
Put $\tV \= \cup_{c \in \CJ} \tV^c$.
Fix $j \in \{0, \ldots, \ell \}$ and note that Part~$3$ of Proposition~\ref{p:BC-UP} implies that $R$ is \BC{$\tau^j\delta_0$}{$\dconst_0\delta$}.
Moreover Lemma~\ref{l:pre construction} applied to $\tdelta = \tau^{-1}\delta_0$, $\teta = 2$, $\delta' = \dconst_0\delta$ and $\hdelta = \tau^j\delta_0$, implies that for each $v \in \CVJ$ there is a simply-connected neighborhood $U_j^v$ of $v$ in $\CC$ such that $\partial U_j^v \subset K(B(\CJ, \tau^{-1}\delta_0))$ and such that
$$
B(v, \tau^j\delta_0) \subset U_j^v \subset B(v, \eta\tau^j\delta_0).
$$
For $c \in \CJ$ let $V_j^c$ be the connected component of $R^{-1}(U_j^{R(c)})$ containing~$c$ and $V_j = \cup_{c \in \CJ} V_j^c$.
Then the assertion follows from Lemma~\ref{l:nice sets with closures}, considering that $\eta\delta_0 < \tau^{-1}\delta_0$.

\partn{2}
Lemma~\ref{l:prenice} applied to $\delta_0$ instead of $\tdelta$, implies that for each $c \in \CJ$ the connected component $\tV^c$ of $\CC \setminus K(B(\CJ, \delta_0))$ containing~$c$ satisfies,
$$
\tB(c, \delta_0) \subset \tV^c \subset \tB(c, 2\delta_0).
$$
For each $j \in \{ 1, \ldots, \ell + 1\}$, Lemma~\ref{l:pre construction} applied to $\tdelta = \delta_0$, $\teta = 2$, $\delta = \tau^{j} \delta_0$, $\delta' = \dconst_0\tau^{\ell + 1} \delta_0$ and $\hdelta = \tau^j\delta_0$, implies that for each $v \in \CVJ$ there is a simply-connected neighborhood $U_j^v$ of $v$ in $\CC$ such that $\partial U_j^v \subset K(B(\CJ, \delta_0))$ and such that
$$
B(v, \tau^j\delta_0) \subset U_j^v \subset B(v, \eta\tau^j\delta_0).
$$
For $c \in \CJ$ let $V_j^c$ be the connected component of $R^{-1}(U_j^{R(c)})$ containing~$c$ and put $V_j = \cup_{c \in \CJ} V_j^c$.
Lemma~\ref{l:nice sets with closures} implies that for each $j = 1, \ldots, \ell + 1$ the set $V_j$ is a nice set and that $(V_1, \ldots, V_{\ell + 1})$ is a nice nest.

Put $V_0 \= \cup_{c \in \CJ} \tV^c$ and assume by induction that for an integer $j > \ell + 1$, nice sets $V_1, \ldots, V_{j - 1}$ satisfying the desired properties have been constructed, and such that moreover for every $k = 1, \ldots, j - 1$ we have $R(\partial V_{k}) \subset K(V_{k  - 1})$.
We will now construct $V_{j}$.
Lemma~\ref{l:pre construction} applied to $\tV_0 = V_{j - \ell - 1}$, $\tdelta = \tau^{j - \ell -1}\delta_0$, $\teta = \eta \le 2$, $\delta = \tau^j \delta_0$, $\delta' = \dconst_0 \tau^j \delta_0$ and $\hdelta = \tau^j\delta_0$, implies that for each $v \in \CVJ$ there is a simply-connected neighborhood $U_j^v$ of $v$ in $\CC$ such that $\partial U_j^v \subset K(V_{j - \ell - 1})$ and such that
$$
B(v, \tau^j\delta_0) \subset U_j^v \subset B(v, \eta\tau^j\delta_0).
$$
For $c \in \CJ$ let $V_j^c$ be the connected component of $R^{-1}(U_j^{R(c)})$ containing~$c$ and put $V_j = \cup_{c \in \CJ} V_j^c$.
By induction hypothesis we have that $R(\partial V_{j - \ell}) \subset K(V_{j - \ell - 1})$ and by construction we have $R(\partial V_{j}) \subset K(V_{j - \ell - 1})$.
So Lemma~\ref{l:nice sets with closures} implies that $(V_{j - \ell}, V_{j})$ is a nice couple.
As by induction by hypothesis $(V_{j - \ell}, \ldots, V_{j - 1})$ is a nice nest, it follows that $(V_{j - \ell}, \ldots, V_j)$ is a nice nest.
Finally note that $K(V_{j - \ell - 1}) \subset K(V_{j - 1})$, so $R(\partial V_j) \subset K(V_{j - 1})$ and the induction hypothesis is satisfied.
\end{proof}

\section{\CL.}\label{s:connecting lemma}
In this section we state and prove a strengthened version of the \CL.
The statement is in~\S\ref{ss:strong connecting lemma} and its proof, that depends on some key area estimates proved in~\S\ref{ss:area estimates}, is in~\S\ref{ss:proof of Theorem D}.

Throughout all this section we fix a rational map $R$ of degree at least~$2$.
\subsection{Strengthened version of the \CL}\label{ss:strong connecting lemma}
Recall that for a constant $\dconst > 0$ we say that $R$ is \BCg{} with constant~$\dconst$, if for every $\delta > 0$ small $R$ is \BC{$\delta$}{$\dconst\delta$}.

For $z \in \CC$ we denote by $\omega(z)$ the $\omega$-limit set of $z$.
\begin{theoalph}[\CL]\label{t:strong connecting lemma}
There are constants $\dconst > 1$ and $K > 1$, only depending on $R$, such that if $R$ is \BCg{} with constant~$\dconst$ and the set
\begin{equation}\label{e:non accumulating}
\{ z \in \CC \mid \omega(z) \cap (\CJ) = \emptyset \}
\end{equation}
has positive Lebesgue measure, then for every sufficiently small $\delta > 0$ the following properties hold.
\begin{enumerate}
\item[1.]
For each $c \in \CJ$ there is $v_c \in \CC$ such that $\dist(v_c, R(c)) < \delta$ and such that for every $n \ge 1$ we have $R^n(v_c) \not\in \tB(\CJ, 2\delta)$.
\item[2.]
For each $c \in \CJ$ let $v_c$ be as in part~$1$ and let $\xi : \CC \to \CC$ be a qc homeomorphism that coincides with the identity outside $\tB(\CJ, \delta)$ and such that for every $c \in \CJ$ we have $\xi(R(c)) = v_c$.
Then there is a rational map $Q$ of the same degree as $R$ and a continuous map $h : \CC \to \CC$ such that the map $\tR \= \xi \circ R : \CC \to \CC$ satisfies
$$
h \circ \tR \equiv Q \circ h
\text{ on } \CC.
$$
\item[3.]
Let $\tR$, $Q$ and $h$ be as in part~$2$.
Then there is a nice set $V = \cup_{c \in \CJ} V^c$ such that for every $c \in \CJ$ we have $\tB(c, \delta) \subset V^c \subset \tB(c, 2\delta)$ and such that there exists a $K$-qc $V$~pseudo-conjugacy between $R$ and $Q$.
\end{enumerate}
\end{theoalph}
By the Fatou-Sullivan classification of connected components of the Fatou set~\cite{CG,Milbook}, it follows that when $J(R) \neq \CC$, the set~\eqref{e:non accumulating} has non empty interior, and hence positive Lebesgue measure.
See also Corollary~\ref{c:expansion away}.
Thus the \CL{} stated in the introduction is a direct consequence of Theorem~\ref{t:strong connecting lemma}.

The following corollary is a direct consequence of Theorem~\ref{t:strong connecting lemma}.
\begin{generic}[Corollary of Theorem~\ref{t:strong connecting lemma}]
Let $R$ be as in Theorem~\ref{t:strong connecting lemma}.
Then for every $\delta > 0$ there is a nice set $V$ for $R$ contained in $\tB(\CJ, \delta)$, such that $R$ is $V$~pseudo-conjugated to a rational map for which no critical point in its Julia set accumulates on a critical point under forward iteration.
\end{generic}
\subsection{Area estimates.}\label{ss:area estimates}
The purpose of this subsection is to prove the following proposition.

For a subset $X$ of $\CC$, we denote by $|X|$ the spherical area of $X$.
\begin{prop}\label{p:area estimates}
Suppose that the set~\eqref{e:non accumulating} has positive Lebesgue measure.
Then for every $\varepsilon \in (0, 1)$ there is $\modul > 0$ such that the following property holds.
Let $( V_j )_{j \ge 1}$ be a nested sequence of nice sets such that for every $j \ge 1$ the sequence $(V_{j}, V_{j + 1}, V_{j + 2}, V_{j + 3})$ is a nice nest and such that for every integer $j \ge 1$ and every $c \in \CJ$ we have $\modulus(V_j^c \setminus \ov{V_{j + 1}^c}) \ge \modul$.
Then for every $c \in \CJ$ and every sufficiently large $j$ we have
$$
\frac{|V_j^c \setminus K(V_{j + 1})|}{|V_j^c|} < \varepsilon.
$$
\end{prop}
The rest of this subsection is dedicated to the proof of this proposition.

Given a nice set $V$ for $R$ and a subset $X$ of $\CC$ of positive Lebesgue measure, put
$$
\xi(X, V) = \frac{|X \setminus R^{-1}(K(V))|}{|X|}
\text{ and }
\psi(X, V) = 1 - \xi(X, V).
$$
Moreover, for a nice set $\hV \= \cup_{c \in \CJ} \hV^c$ for $R$, such that $(\hV, V)$ is a nice couple, put
$$
\Xi(\hV, V) = \max_{c \in \CJ} \frac{|\hV^c \setminus K(V)|}{|\hV^c|}
\text{ and }
\Psi(\hV, V) = 1 - \Xi(\hV, V).
$$
Finally, given a constant $D > 1$ put,
$$
\tXi(\hV, V) \= \min \{ D^2 \Xi (\hV, V), 1 - D^{-2} \Psi (\hV, V) \},
\text{ and }
\tPsi(\hV, V) = 1 - \tXi(\hV, V).
$$
Note that we have $\tXi(\cdot,\cdot), \tPsi(\cdot,\cdot) \in [0, 1]$, $D^{-2}\tXi(\cdot,\cdot) \le \Xi(\cdot,\cdot) \le \tXi(\cdot,\cdot)$ and $D^{-2}\tPsi(\cdot,\cdot) \le \Psi(\cdot,\cdot) \le \tPsi(\cdot,\cdot)$.
Let $m$ be a positive integer and let $\hW$ be a connected component of $R^{-m}(\hV)$ such that $R^m$ is univalent on $\hW$ and such that the distortion of $R^m$ on $\hW$ is at most~$D$.
Then we have
$$
\xi(\hW, V) \le \tXi(\hV, V).
$$
\begin{lemm}\label{l:distribution}
Let $D > 1$ and let $(\hV, V)$ be a nice couple for $R$.
Moreover, let $X$ be a subset of $\CC$ of positive Lebesgue measure satisfying the following property: For every connected component $\hW$ of $\CC \setminus R^{-1}(K(\hV))$ that intersects $X$ we have that $\hW \subset X$ and that $R^{m_{\hW}}$ is univalent on $\hW$ with distortion bounded by $D$.
Then,
$$
\xi(X, V) \le \xi(X, \hV) \cdot \tXi(\hV, V).
$$
\end{lemm}
\begin{proof}
Denote by $\sC$ (resp. $\hsC$) the collection of connected components of $\CC \setminus R^{-1}(K(V))$ (resp. $\CC \setminus R^{-1}(K(\hV))$).
Then
$$
\sum_{W \in \sC, W \subset \hW} |W| = |\hW \setminus K(V)| = |\hW| \xi(\hW, V) \le |\hW| \cdot \tXi (\hV, V).
$$
Thus
\begin{multline*}
|X| \xi(X, V) = |X \setminus K(V)| = \sum_{W \in \sC, W \subset X} |W| = \sum_{\hW \in \hsC, \hW \subset X} \sum_{W \in \sC, W \subset \hW} |W|
\le \\ \le
\left( \sum_{\hW \in \hsC, \hW \subset X} |\hW| \right) \tXi(V, \hV) \le |X| \xi(X, \hV) \cdot \tXi(\hV, V).
\end{multline*}
\end{proof}

Let $\modul > 0$ and let $D(\modul) > 1$ be the constant defined in~\S\ref{ss:nice nests}.
Then, if $(\hhV, \hV, V)$ is a $\modul$-separated nice triple for $R$ and $W$ is a \cc{} of $(\hV \setminus \ov{V}) \setminus R^{-1}(K(\hV))$ (resp. $(\hV \setminus \ov{V}) \setminus R^{-1}(K(\hhV))$), then the map $R^{m_W}$ is univalent on $W$ and its distortion on $W$ is at most $D(\modul)$.
\begin{lemm}\label{l:inductive estimate}
Let $\modul > 0$, let $(\hhV, \hV, V)$ be a $\modul$-separated nice triple and put
$$
D \= D(\modul)
\text{ and }
\sigma \= \max_{c \in \CJ} \frac{|\hV^c|}{|\hV^c \setminus \ov{V^c}|}.
$$
Then we have
$$
\tPsi(\hV, V) \ge \frac{\tPsi(\hhV, \hV)}{D^2 \sigma - 1 + \tPsi(\hhV, \hV)}.
$$
\end{lemm}
\begin{proof}
Given $c \in \CJ$, Lemma~\ref{l:distribution} applied to the nice couple $(\hV, V)$ and to $X = \hV^c \setminus \ov{V^c}$ gives,
$$
\xi(\hV^c \setminus \ov{V^c}, V)
\le
\xi(\hV^c \setminus \ov{V^c}, \hV) \cdot \tXi(\hV, V).
$$
Similarly, Lemma~\ref{l:distribution} applied to the nice couple $(\hhV, \hV)$ and to $X = \hV^c \setminus \ov{V^c}$ gives
$$
\xi(\hV^c \setminus \ov{V^c}, \hV)
\le
\xi(\hV^c \setminus \ov{V^c}, \hhV) \cdot \tXi(\hhV, \hV)
\le
\tXi(\hhV, \hV).
$$
Combining these two inequalities we get
$$
\xi(\hV^c \setminus \ov{V^c}, V)
\le
\tXi(\hhV, \hV) \cdot \tXi(\hV, V),
$$
and
$$
\psi(\hV^c \setminus \ov{V^c}, V)
\ge
\tPsi(\hhV, \hV) + \tPsi(\hV, V) - \tPsi(\hhV, \hV) \cdot \tPsi(\hV, V).
$$
As $V$ is disjoint from $K(V)$ we have,
$$
\psi(\hV^c \setminus \ov{V^c}, V)
=
\frac{|\hV^c \cap K(V)|}{|\hV^c \setminus \ov{V^c}|}
=
\frac{|\hV^c|}{|\hV^c \setminus \ov{V^c}|} \cdot \frac{|\hV^c \cap K(V)|}{|\hV^c|}.
$$
Therefore,
\begin{multline*}
\sigma D^2 \tPsi(\hV, V) \ge \sigma \Psi(\hV, V)
\ge
\min_{c \in \CJ} \psi(\hV^c \setminus \ov{V^c}, V)
\ge \\ \ge
\tPsi(\hhV, \hV) + \tPsi(\hV, V) - \tPsi(\hhV, \hV) \cdot \tPsi(\hV, V),
\end{multline*}
which easily implies the desired inequality.
\end{proof}

\begin{proof}[Proof of Proposition~\ref{p:area estimates}]
Let $\modul > 0$ be sufficiently large so that $D(\modul) < (1 + \varepsilon)^{\frac{1}{3}}$ and so that every annulus in $\CC$ of modulus at least $\modul$ intersects the Julia set of $R$ (see~\cite{MR}).
Taking $\modul$ larger if necessary we have that, if $( V_j )_{j \ge 1}$ is as in the statement of the proposition, then
$$
\sigma \= \sup_{j \ge 1, c \in \CJ} \frac{|V_j^c \setminus \ov{V_{j + 1}^c}|}{|V_j^c|} < (1 + \varepsilon)^{\frac{1}{3}}.
$$
Note in particular that $D(\modul)^2\sigma - 1 < \varepsilon$ and that for every $j \ge 1$ and $c \in \CJ$ the annulus $V_j^c \setminus \ov{V_{j + 1}^c}$ intersects the Julia set of $R$.

\partn{1}
We will prove that for every sufficiently large $j$ and every $c \in \CJ$ the set $V_j^c \setminus K(V_{j + 1})$ has positive Lebesgue measure.

\partn{1.1}
Let us prove by induction that for every $c \in \CJ$ and $j \ge 1$ there is $c' \in \CJ$, an integer $n \ge 0$ and a \cc{} $W$ of $R^{-n}(V_1^c)$ that is contained in $V_j^c \setminus \ov{V_{j + 1}^c}$ and such that for every $m = 1, \ldots, n - 1$ we have $R^m(W) \cap V_j = \emptyset$.
When $j = 1$ just take $c' = c$, $n = 0$ and $W = V_1^c$.
Suppose by induction that this property holds for some $j \ge 1$ and let $c \in \CJ$.
Since the annulus $V_{j + 1}^c \setminus \ov{V_{j + 2}^c}$ intersects the Julia set of $R$, it follows by the eventually onto property of Julia sets that the set $(V_{j + 1}^c \setminus \ov{V_{j + 2}^c}) \setminus R^{-1}(K(V_{j}))$ is non empty.
Let $W'$ be a \cc{} of this set.
Then $R^{m_{W'}}$ maps $W'$ univalently to $V_j^{c(W')}$ and for every $m = 1, \ldots, m_{W'} - 1$ the set $R^m(W')$ is disjoint from $V_{j + 1}$.
Let $W$ be given by the induction hypothesis for $j$ and $c = c(W')$, so that $W \subset V_j^c \setminus \ov{V_{j + 1}^c}$.
Then the set $R^{m_{W'}}|_{W'}^{-1}(W)$ satisfies the desired properties.

\partn{1.2}
Note that the set~\eqref{e:non accumulating} is equal to the union of $\cup_{j \ge 1} K(V_j)$ and of the iterated preimages of the non recurrent critical points of $R$.
So for sufficiently large $j$ the set $K(V_j)$ has positive Lebesgue measure. 
It follows that for every repelling periodic point $p$ of $R$ there is a positive integer $j(p)$ such that $p \in K(V_{j(p)})$ and such that the Lebesgue density of $K(V_{j(p)})$ is positive at $p$.
As for each $c \in \CJ$ the annulus $V_1^c \setminus \ov{V_2^c}$ intersects the Julia set of $R$, there is a repelling periodic point $p_c$ of $R$ contained in this set.
Put $j_0 \= \max_{c \in \CJ}j(p_c)$ and let $j \ge j_0$ be given.
By part~$1.1$ it follows that for a given $c \in \CJ$ there is a point $p$ in $V_j^c \setminus \ov{V_{j + 1}^c}$, an integer $n \ge 0$ and $c' \in \CJ$ such that $R^n(p) = p_{c'}$ and such that for every $m = 0, \ldots, n - 1$ we have $R^m(p) \not\in V_{j}$.
As $j \ge j_0$ and $p_c \in K(V_{j(p_c)}) \subset K(V_{j_0})$, it follows that $p \in K(V_j)$ and that the Lebesgue density of $K(V_{j + 1})$ is positive at~$p$.
This proves that the set $(V_j^c \setminus \ov{V_{j + 1}^c}) \cap K(V_{j + 1})$ has positive Lebesgue measure.

\partn{2}
For $j \ge 1$ put $\Psi_j = \Psi(V_j, V_{j + 1})$ and $\tPsi_j \= \tPsi(V_j, V_{j + 1})$.
We need to prove that $\limsup_{j \to \infty} \Psi_j > 1 - \varepsilon$.
As $\tPsi_j \le \Psi_j$, we just need to prove that $\limsup_{j \to \infty} \tPsi_j > 1 - \varepsilon$.

Note that by hypothesis for each $j \ge 2$ the nice triple $(V_{j + 2}, V_{j + 1}, V_{j})$ is $\modul$-separated.
Thus Lemma~\ref{l:inductive estimate} implies that for every $j \ge 2$ we have
\begin{equation}\label{e:inductive estimate}
\tPsi_{j + 1} \ge \tPsi_{j}(D(\modul)^2 \sigma - 1 + \tPsi_{j})^{-1}.
\end{equation}
By our choice of $\modul$ we have $D(\modul)^2 \sigma - 1  < \varepsilon$.
Since by part~$1$ for every sufficiently large $j$ we have $\tPsi_j \ge D^{-2} \Psi_j > 0$, the inequality \eqref{e:inductive estimate} easily implies that $\liminf_{j \to \infty} \tPsi_j > 1 - \varepsilon$.
\end{proof}
\subsection{Proof of Theorem~\ref{t:strong connecting lemma}}\label{ss:proof of Theorem D}
Note that the usual \KDT{} implies that for every constant $\modul > 0$ there is $\tD(\modul) > 1$ such that the following property holds.
Let $U$ be as subset of $\CC$ biholomorphic to $\D$ and such that $\diam(U) < \tfrac{1}{2} \diam(\CC)$.
Moreover, let $\cK$ be a compact subset of $U$ such that $U \setminus \cK$ is an annulus of modulus at least $\modul$.
Then for every map $\varphi : U \to \C$ that is univalent on $\cK$, the distortion of $\varphi$ on $\cK$ is bounded by~$\tD(\modul)$.
\begin{lemm}\label{l:absolute area}
Let $(\hhV, \hV, V)$ be a nice triple for $R$ and put
$$
\modul \= \min_{c \in \CJ} \modulus(\hhV^c \setminus \ov{\hV^c}).
$$
Then for every $c_0 \in \CJ$ we have,
$$
\left\| \hhV^{c_0} \setminus K(V) \right\|_{\hhV^{c_0}} \le |\D| \tD(\modul)^2 \max_{c \in \CJ} \frac{|\hV^c \setminus K(V)|}{|\hV^c|}.
$$
\end{lemm}
\begin{proof}
Given $c_0 \in \CJ$ consider a biholomorphism $\varphi : \hhV^{c_0} \to \D$.
For each connected component $W$ of $\hhV^{c_0} \setminus K(\hV)$ there is $m_W \ge 0$ and $c(W) \in \CJ$ such that $R^{m_W} : W \to \hV^{c(W)}$ is a biholomorphism whose inverse $\phi_W$ extends univalently to $\hhV^{c(W)}$.
So the distortion of $\varphi \circ \phi_W$ on $\hhV^{c_0}$ is bounded by $\tD(\modul)$.
Therefore,
$$
\frac{|\varphi(W \setminus K(V))|}{|\varphi(W)|} \le 
\tD(\modul)^2 \frac{|\hV^{c(W)} \setminus K(V)|}{|\hV^{c(W)}|}.
$$
Since the connected components of $\hhV^{c_0} \setminus K(\hV)$ cover $\hhV^{c_0} \setminus K(V)$, the assertion of the lemma follows.
\end{proof}
\begin{lemm}\label{l:verifying rigidity}
Let $\tau \in (0, 1)$ and $\eta \in (1, \tau^{-1})$ be given.
Then there are $K > 1$ and a positive integer $\ell$, only depending on $\tau$ and $R$, such that the following property holds.
Let $\delta_0 > 0$ be small and let $( V_j )_{j \ge 1}$ be a  nested sequence of nice sets such that for every $j \ge 1$ the sequence $(V_{j}, \ldots, V_{j + \ell})$ is a nice nest, and such that for every $c \in \CJ$ we have
$$
\tB(c, \tau^j \delta_0) \subset V_j^c \subset \tB(c, \eta\tau^j \delta_0).
$$
Then for every sufficiently large $j$ the nice couple $(V_{j}, V_{j + 1})$ is $K$-rigid and for each $c \in \CJ$ the set $V_j^c \cap K(V_{j + 1}^c)$ has positive Lebesgue measure.
\end{lemm}
\begin{proof}
Let $\modul_0 > 0$ and $\varepsilon_0 > 0$ be given by Rigidity with $K = 2$.
Moreover let $\modul_1 > 0$ be the constant $\modul$ given by Proposition~\ref{p:area estimates} for $\varepsilon \= (|\D|\tD(\modul_0)^2)^{-1}\varepsilon_0$.
Let $N$ be a positive integer such that if $( V_j )_{j \ge 1}$ is as in the statement of the lemma, then for every integer $j \ge 1$ and every $c \in \CJ$ we have
$$
\modulus(V_j^c \setminus \ov{V_{j + N}^c}) \ge \max \{ \modul_0, \modul_1 \}.
$$
Put $\ell = 3N$ and notice that for every $j_0 \ge 1$ the sequence of nice nests $( V_{j_0 + kN} )_{k \ge 0}$ satisfies the hypothesis of Proposition~\ref{p:area estimates}.
So for every sufficiently large $j$ and every $c \in \CJ$ we have
\begin{equation}\label{e:first estimate}
|V_j^c \setminus K(V_{j + N})| < \varepsilon |V_j^c|.
\end{equation}
In particular, for a given $c \in \CJ$, Lemma~\ref{l:distribution} applied to the nice couple $(V_{j}, V_{j + N})$ and to $X = V^c_j \setminus \ov{V^c_{j + 1}}$ implies that the set $V_j^c \cap K(V_{j + 1})$ has positive Lebesgue measure.

On the other hand, \eqref{e:first estimate} and Lemma~\ref{l:absolute area} applied to the nice triple $(V_{j}, V_{j + N}, V_{j + 2N})$ imply that,
$$
\| V_j^c \setminus K(V_{j + 2N}) \|_{V_j^c} < \varepsilon_0.
$$
So Rigidity applied to the nice triple $(V_{j}, V_{j + 2N}, V_{j + 3N})$ implies that the nice couple $(V_{j}, V_{j + 2N})$ is $2$-rigid.
Finally part~$2$ of Proposition~\ref{p:rigid nice couples} implies that for some constant $K' > 1$ only depending on $\tau$ and $R$, the nice couple $(V_{j}, V_{j + 1})$ is $2K'$-rigid.
\end{proof}
\begin{proof}[Proof of Theorem~\ref{t:strong connecting lemma}]
Put $\tau = 2^{\frac{1}{4}}$, choose $\eta \in (0, \tau^{-1})$ and let $\ell_0$ be the integer $\ell$ given by Lemma~\ref{l:verifying rigidity} for this choice of $\tau$.
Moreover, let $\dconst = \dconst_0 > 0$ be the constant as in Proposition~\ref{p:nice nest}, for this choice of $\tau$, $\eta$ and for $\ell \= \max \{ 4, \ell_0 \}$.

Suppose that $R$ is \BCg{} with constant~$\dconst$.
Then, by part~$2$ Proposition~\ref{p:nice nest} there is $\delta_0 > 0$ and a sequence of nice sets $( V_j )_{j \ge 0}$ satisfying the hypothesis of Lemma~\ref{l:verifying rigidity}.
So for every $\delta > 0$ small there is $j \ge 1$ such that for every $c \in \CJ$ we have
$$
\tB(c, \delta) \subset V_{j + 2}^c \subset 
V_j^c \subset \tB(c, 2\delta),
$$
and that the nice triple $(V_{j}, V_{j + 1}, V_{j + 2})$ verifies the hypothesis of Theorem~\ref{t:Thurston's algorithm}.
Parts~$2$ and~$3$ are then direct consequences of Theorem~\ref{t:Thurston's algorithm}.
It reminds to prove part~$1$.

To prove part~$1$, fix $c \in \CJ$.
Note that we can choose $j$ such that in addition we have $\tB(\CJ, 2\delta) \subset V_{j - 2}$.
It follows that for every point $z \in K(V_{j - 2})$ and every positive integer $n$, we have $R^n(z) \not\in \tB(\CJ, 2\delta)$.
Thus we just need to show that there is a point $v_c$ in $K(V_{j - 2})$ whose distance to $R(c)$ is less than $\delta$. 
If $R(c) \in K(V_{j - 2})$ then we can take $v_c = R(c)$, so we assume that $R(c) \not \in K(V_{j - 2})$.
Denote by $W$ the \cc{} of $\CC \setminus K(V_{j - 2})$ containing $R(c)$.
As $\ell \ge 4$ it follows that $(V_{j - 2}, V_{j + 2})$ is a nice couple and therefore that the closure of the \cc{} of $R^{-1}(W)$ containing $c$ is contained in $V_{j + 2}$.
Thus $\ov{W} \subset R(V_{j + 2}^c) \subset B(R(c), \delta)$.
So each point in $\partial W$ belongs to $K(V_{j - 2})$ and its distance to $R(c)$ is less than $\delta$.
This completes the proof of part~$1$ and of the theorem.
\end{proof}

\section{Expanding properties of \BCg{} rational maps.}\label{s:global}
In this section we show that a \BCg{} rational map has no Siegel discs nor Herman rings.
This and other results are proved in~\S\ref{ss:expansion away}, using a theorem of Ma\~n\'e~\cite{Mane, SL} and the Corollary of Theorem~D.
Moreover we prove Theorem~\ref{t:measure} in~\S\ref{ss:proof of Theorem B}.
\subsection{Expansion away of critical points.}\label{ss:expansion away}
The purpose of this section is to prove the following result, which is a rather easy consequence of Theorem~\ref{t:strong connecting lemma}.
\begin{prop}\label{p:expansion away}
Let $R$ be a rational map of degree at least~$2$, let $\dconst > 0$ be the constant given by Theorem~\ref{t:strong connecting lemma} and suppose that $R$ is \BCg{} with constant~$\dconst$.
If the set~\eqref{e:non accumulating} has positive Lebesgue measure, then the following properties hold.
\begin{enumerate}
\item[1.]
$R$ does not have irrationally indifferent cycles nor Herman rings.
\item[2.]
The set $\{ z \in J(R) \mid \omega(z) \cap \Crit(R) = \emptyset \}$ has zero Lebesgue measure.
\item[3.]
If $K \subset J(R)$ is a compact and forward invariant set disjoint from $\Crit(R)$, then $R$ is expansive on $K$.
If moreover $K$ contains no parabolic periodic point of $R$, then $R$ is uniformly expanding on $K$.
\end{enumerate}
\end{prop}
By the Fatou-Sullivan classification of connected component of the Fatou set~\cite{CG,Milbook} it follows that when $J(R) \neq \CC$, the set~\eqref{e:non accumulating} has non empty interior, and hence positive Lebesgue measure.
So the following corollary is a direct consequence of part~$2$ of the proposition.
\begin{coro}\label{c:expansion away}
Let $\dconst > 0$ be given by Theorem~\ref{t:strong connecting lemma} and suppose that $R$ is \BCg{} with constant~$\dconst$.
Then the set~\eqref{e:non accumulating} has positive Lebesgue measure if and only if $J(R) \neq \CC$.
\end{coro}
The following corollary is a direct consequence of the previous corollary and of Proposition~\ref{p:expansion away}.
\begin{coro}\label{c:expansion away 2}
Let $\dconst > 0$ be given by Theorem~\ref{t:strong connecting lemma} and suppose that $R$ is \BCg{} with constant~$\dconst$.
Then the following properties hold.
\begin{enumerate}
\item[1.] $R$ has no Siegel disks nor Herman Rings.
If moreover $J(R) \neq \CC$, then $R$ has no Cremer periodic points.
\item[2.]
If $J(R) \neq \CC$ and $K \subset J(R)$ is a compact and forward invariant set disjoint from $\Crit(R)$, then $R$ is expansive on $K$.
If moreover $K$ contains no parabolic periodic point of $R$, then $R$ is uniformly expanding on $K$.
\end{enumerate}
\end{coro}
The proof of Proposition~\ref{p:expansion away} is at the end of this subsection.
It is based on the following general lemma.
\begin{lemm}\label{l:pseudo-conjugacy}
Let $R$ and $Q$ be rational maps of degree at least~$2$, let $V$ be a nice set for $R$ and let $\chi$ be a $V$~pseudo-conjugacy between $R$ and $Q$.
Then the following properties hold.
\begin{enumerate}
\item[1.]
If $\cO$ is a periodic orbit of $R$ that is disjoint from $V$, then $\chi(\cO)$ is a periodic orbit of $Q$ of the same nature as $\cO$ (repelling, attracting, parabolic, Siegel or Cremer.)
\item[2.]
A Siegel disk or a Herman ring of $R$ (resp. $Q$) that intersects $K(V)$ (resp. $\chi(K(V))$), is contained in $K(V)$ (resp. $\chi(K(V))$.)
In particular the image by $\chi$ of a Siegel disk or a Herman ring of $R$ intersecting $K(V)$, is a Siegel disk or a Herman ring of~$Q$.
\item[3.]
A point in $K(V)$ belongs to the Fatou set of $R$ if and only if its image by~$\chi$ belongs to the Fatou set of $Q$.
\item[4.]
Let $K \subset \CC$ be a compact set disjoint from $\ov{V}$ that is forward invariant by~$R$.
If $Q$ is uniformly expanding on $\chi(K)$, then~$R$ is uniformly expanding on~$K$.
\end{enumerate}
\end{lemm}
\begin{proof}
\

\partn{1}
As $V$ is a nice set of $R$, it follows that $\partial V$ does not contain periodic points of~$R$ (cf. Lemma~\ref{l:nice sets forward}.)
Thus the periodic orbit $\cO$ is contained in the open set $\CC \setminus \ov{V}$, and $\chi$ conjugates $R$ on a neighborhood of $\cO$ to $Q$ on a neighborhood of $\chi(\cO)$.
The dynamics on a neighborhood of a repelling, attracting, parabolic or Siegel periodic orbit can be easily characterized in topological terms, so the periodic orbits $\cO$ and $\chi(\cO)$ have the same nature.

\partn{2}
We will prove the assertion for $R$; the proof of the assertion for $Q$ is analogous.
Note that it is enough to prove that every Siegel disk and every Herman ring of $R$ is disjoint from $\partial V$.

Suppose by contradiction that there is a Siegel disk or a Herman ring $\cR$ of $R$ that contains a point $z_0$ in~$\partial V$.
Let $n$ be the period of $\cR$ and for $z \in \CC$ denote by $\omega_n(z)$ the $\omega$-limit set of $z$ under $R^n$.
Then for each $z \in \cR$ set $\omega_n(z)$ is an analytic Jordan curve, or we have $R^n(z) = z$ and $\omega_n(z) = \{ z \}$.
In both cases $R^n(\omega_n(z)) = \omega_n(z)$ and $R^n$ is transitive on $\omega_n(z)$.
Moreover note that $\omega_n(z)$ depends continuously with $z$ in the Hausdorff topology.

As $V$ is a nice set, the forward orbit of each point in $\partial V$ is disjoint from~$\ov{V}$.
It follows that for every $z \in \cR$ the set $\omega_n(z)$ is either contained in~$V$ or disjoint from~$V$.
Thus $\omega_n(z_0)$ is disjoint from~$V$. 
As $\cR$ is open, the set $\omega_n(z_0)$ can be approximated in the Hausdorff topology by sets of the form $\omega_n(z)$, with $z \in \cR \cap V$.
As for such $z$ we have $\omega_n(z) \subset \cR \cap V$, we conclude that $\omega_n(z_0) \subset \ov{V}$.
Since $\omega_n(z_0)$ is disjoint from $V$, we then have $\omega_n(z_0) \subset \partial V$ and we get a contradiction with the fact that $V$ is a nice set.

\partn{3}
We will just prove the direct implication; the proof of the reverse implication being analogous.

Let $p \in \CC$ be in the Fatou set of $R$.
By the Fatou-Sullivan classification of connected components of the Fatou set, we have the following two cases.

\partn{Case 1}
\textsf{\textit{For some integer $n \ge 1$ the point $R^n(p)$ belongs to a Siegel disk or a Herman ring of $R$.}}
Then we have $R^n(p) \in K(V)$ and part~$2$ implies that $\chi(R^n(p)) = Q^n(\chi(p))$ belongs to a Siegel disk or a Herman ring of $Q$.
In particular $\chi(p)$ belongs to the Fatou set of $Q$.

\partn{Case 2}
\textsf{\textit{The forward orbit of $p$ is asymptotic to an attracting or parabolic periodic orbit $\cO$ of $R$.}}
Then we have $\cO \subset K(V)$ and by part~$1$ it follows that $\chi(\cO)$ is an attracting or parabolic periodic orbit of $Q$.
As $\chi$ is a $V$~pseudo-conjugacy, it follows that $\chi(p)$ is asymptotic to $\chi(\cO)$ and that $\chi(p)$ belongs to the Fatou set of $Q$.

\partn{4}
As $Q$ is uniformly expanding on $\chi(K)$, it follows that $Q$ does not have critical points in $\chi(K)$.
As $\chi$ conjugates $R$ and $Q$ on $\CC \setminus V$ and $K \subset \CC \setminus \ov{V}$, it follows that $R$ does not have critical points in $K$.
So there is $\varepsilon_0 > 0$ such that for all $w \in K$ the distortion of $R$ on $B(w, \varepsilon_0)$ is bounded by some constant $D > 1$, close to~1.

Since $Q$ is uniformly expanding on $\chi(K)$, there are $\varepsilon_1 > 0$, $C_1 > 0$ and $\lambda_1 > 1$ such that for all $z \in \chi(K)$ the pull-back $W$ of $B(Q^n(z), \varepsilon_1)$ by $Q^n$ to $z$ is univalent and $\diam(W) \le C_1 \lambda_1^{-n}$.
We assume $\varepsilon_1 > 0$ small enough so that $\diam(\chi^{-1}(W)) < \varepsilon_0$.
By the H\"older property of qc homeomorphisms there are $\alpha \in (0,1)$ and $C_0 > 0$ be such that for all $z_0, z_1 \in \CC$ we have
$$
\dist(\chi^{-1}(z_0), \chi^{-1}(z_1)) \le C_0 \dist(z_0, z_1)^\alpha,
$$
see Appendix~\ref{a:qc}.
Suppose that $\varepsilon_1 \ll \diam(\chi(K(V)))$ and fix $z \in \chi(K)$.
Since $K(V)$ is connected (Lemma~\ref{l:nice sets forward}), for $n$ big there is
$$
w_n \in \partial B(Q^n(z), \varepsilon_1 / 2) \cap \chi(K(V)),
$$
so that $\dist(\chi^{-1}(Q^n(z)), \chi^{-1}(w_n)) \sim 1$.

Consider the preimage $y_n$ of $w_n$ by $Q^n$ in the connected component of $Q^{-n}(B(Q^n(z),\varepsilon_1))$ that contains $z$.
So $\dist(z, y_n) \le C_1 \lambda_1^{-n}$ and therefore $\dist(\chi(z), \chi(y_n)) \le C_0C_1^\alpha \lambda_1^{-\alpha n}$.
So there is a constant $C > 0$ such that $|(R^n)'(z)| \ge C(\lambda_1^\alpha D^{-1})^n$.
Taking $\varepsilon_0$ smaller if necessary we may assume that $\lambda_1^\alpha D^{-1} > 1$.
\end{proof}

\begin{proof}[Proof of Proposition~\ref{p:expansion away}]
\

\partn{1}
Assume by contradiction that $R$ has an irrationally indifferent periodic point $p$ or a Herman ring $H$.
In the latter case choose $p \in H$.
Let $\delta > 0$ be sufficiently small, so that $\tB(\CJ, \delta)$ is disjoint from the forward orbit of $p$. 
By the Corollary of Theorem~\ref{t:strong connecting lemma} there is a nice set $V \subset \tB(\CJ, \delta)$, a rational map $Q$ without recurrent critical points in its Julia set and a $V$~pseudo-conjugacy $\chi$ between $R$ and~$Q$.
Parts~$1$ and~$2$ of Lemma~\ref{l:pseudo-conjugacy} imply that $\chi(p)$ is an irrationally indifferent periodic point of $Q$ or that $\chi(p)$ is contained in a Siegel disk or a Herman ring of $Q$.
Since $Q$ does not have irrationally indifferent periodic points nor Herman rings, we get a contradiction.

\partn{2}
Is enough to prove that for every $\delta > 0$ the set
$$
J(R) \cap K(\tB(\CJ, \delta))
$$
has zero Lebesgue measure.
Let $V \subset \tB(\CJ, \delta)$, $Q$ and $\chi$ be given by the Corollary of Theorem~\ref{t:strong connecting lemma}.
By part~$3$ of Lemma~\ref{l:pseudo-conjugacy} it follows that the image of $J(R) \cap K(V)$ by $\chi$ is contained in the Julia set of $Q$.
As $Q$ has no recurrent critical points, it follows that the image of this set by $\chi$ has zero Lebesgue measure.
As qc homeomorphisms are absolutely continuous with respect to the Lebesgue measure, it follows that $J(R) \cap K(V)$ has zero Lebesgue measure.

\partn{3}
Let $\delta > 0$ be sufficiently small so that the closure of $\tB(\CJ, \delta)$ is disjoint from $K$.
Let $V \subset \tB(\CJ, \delta)$, $Q$ and $\chi$ given by the Corollary of Theorem~\ref{t:strong connecting lemma}.
Then $K \subset K(V)$ and the compact set $\hK = \chi(K)$ is forward invariant by $Q$ and disjoint from the critical points of $Q$.
Moreover, part~$3$ of Lemma~\ref{l:pseudo-conjugacy} implies that $\hK$ is contained in the Julia set of $Q$.

As $Q$ is expansive on $\hK$ (cf. \cite{DU}), it follows that $R$ is expansive on $K$.
If $K$ does not contain parabolic periodic points of $R$, then $\hK$ does not contain parabolic periodic points of~$Q$ (part~$1$ of Lemma~\ref{l:pseudo-conjugacy}.)
In this case $Q$ is uniformly expanding on~$\hK$, so part~$4$ of Lemma~\ref{l:pseudo-conjugacy} implies that $R$ is uniformly expanding on~$K$.
\end{proof}
\subsection{Proof of Theorem~\ref{t:measure}.}\label{ss:proof of Theorem B}
Let $R$ be a rational map of degree at least~$2$.
Let $\dconst > 0$ be the constant given by Theorem~\ref{t:strong connecting lemma} and let $\modul > 0$ be the constant given by Proposition~\ref{p:area estimates} for $\varepsilon = \tfrac{1}{4}$.
Taking $\modul$ larger if necessary we assume that $D(\modul)^2 < 2$.
Let $\tau \in (0, 1)$ be sufficiently small and $\eta \in (1, \tau^{-1})$ sufficiently close to~$1$, so that for every $c \in \CJ$ and every small $\delta > 0$ we have
$$
\modulus \left( \tB(c, \delta) \setminus \ov{\tB(c, \eta\tau\delta)} \right)
\ge \modul.
$$
Let $\dconst_0$ be a constant satisfying the hypothesis of Proposition~\ref{p:nice nest} for these choices of $\tau$ and $\eta$ and for $\ell = 3$.
Suppose that $R$ is \BCg{} with constant $\max \{ \dconst, \dconst_0 \}$.

Then part~$2$ of Theorem~\ref{t:measure} is a direct consequence of Corollary~\ref{c:expansion away}.

To prove part~$1$ of Theorem~\ref{t:measure}, notice first that by part~$2$ of Proposition~\ref{p:expansion away} the set of points in $\CC$ that do not accumulate on a critical point under forward iteration has zero Lebesgue measure.
So we just have to prove that the set of points in $J(R)$ that accumulate on a critical point under forward iteration has zero Lebesgue measure.
Let $\delta_0 > 0$ and $( V_j )_{j \ge 1}$ be the sequence of nice nests given by part~$2$ of Proposition~\ref{p:nice nest}.
Clearly the sets $\CC \setminus K(V_j)$ are decreasing with $j$ and their intersection contains the set of points that accumulates on a critical point under forward iteration.
Thus we just have to prove that the Lebesgue measure of $\CC \setminus K(V_j)$ converges to~$0$ as $j \to \infty$.

It follows by Proposition~\ref{p:area estimates} that for every $c \in \CJ$ and every sufficiently large $j$ we have,
\begin{equation}\label{e:measure}
|V_j^c \setminus K(V_{j + 1})| < \tfrac{1}{4} |V_j^c|.
\end{equation}
Moreover, an easy distortion estimate yields that for every $j \ge 1$ we have,
$$
|\CC \setminus K(V_{j + 1})|
\le
D(\modul)^2 \left(
\max_c \frac{|V^c_j \setminus K(V_{j + 1})|}{|V_j^c|} \right)
|\CC \setminus K(V_j)|.
$$
Since $D(\modul)^2 < 2$, inequality~\eqref{e:measure} implies that for every sufficiently large $j$ we have $|\CC \setminus K(V_{j + 1})| < \tfrac{1}{2} |\CC \setminus K(V_j)|$.
Thus $|\CC \setminus K(V_j)| \to 0$ as $j \to \infty$, and the proof of the theorem is complete.

\appendix
\section{Summability and Collet-Eckmann conditions.}\label{a:summ and CE}
This appendix is dedicated to the proof of Theorem~\ref{t:summ and CE}.
We recall the shrinking neighborhood technique in~\S\ref{ss:shrinking neighborhoods} and, after some preliminary lemmas in~\S\ref{ss:lemmas for summ and CE}, we prove Theorem~\ref{t:summ and CE} in~\S\ref{ss:proof of summ and CE}.

We fix throughout this appendix a rational map $R$ of degree at least~$2$.
\subsection{Shrinking neighborhoods.}\label{ss:shrinking neighborhoods}
The proof of Theorem~\ref{t:summ and CE} is based on a technique to control the distortion of backward iterates, that was introduced by F.~Przytycki in~\cite{Pr1}, see also~\cite{GSmCE,GSm,PrUparabolic,BvS}.
This technique is called \textsf{\textit{shrinking neighborhoods}} and it is described below.
We will use the shrinking neighborhoods only in Lemma~\ref{l:summability0}.

Fix a sequence of numbers $( d_j )_{j \ge 0}$ in $(0, 1)$, such that $\prod_{j \ge 0}(1 - d_j) = \tfrac{1}{2}$.
Put $D_0 = 1$ and for $j > 0$ put $D_j = \prod_{0 \le i \le j -1} (1 - d_i)$.

Consider a critical point $c \in \CJ$, $n \ge 1$ and $\xi \in R^{-n}(c)$.
Given $r \in (0, r_K)$ (where $r_K > 0$ is as in~\S\ref{ss:distortion}) and $j \in \{0, \ldots, n \}$, let $U_j$ (resp. $U_j'$) be the connected component of $R^{-j}(\tB(c, D_{j + 1}r))$ (resp. $R^{-j}(\tB(c, D_j r))$) that contains $R^{n - j}(\xi)$.
Note that $U_j \subset U_j'$.

If for every $j = 1, \ldots, n$ the set $U_j'$ is disjoint from $\Crit(R)$, then Lemma~\ref{l:shrinking neighborhoods} below implies that for every $v \in U_n$ we have
\begin{eqnarray}\label{e:1}
\dist(\xi, v) |(R^n)'(v)| \le K d_n^{-1} \dist(c, R^n(v)).
\end{eqnarray}
Here the constant $K$ only depends on $\mu_{\max}$ and on the constant $r_K > 1$ as in~\S\ref{ss:distortion}.

The following is Lemma~2.1 of~\cite{GSm}.
\begin{lemm}\label{l:shrinking neighborhoods}
Let $U \subset \C$ be conformally equivalent to the unit disc~$\D$ and consider a biholomorphism $f : U \to \D$.
Setting $\xi = f^{-1}(0) \in U$, for every $v \in U$ we have
$$
|\xi - v| \cdot |f'(v)| \le 2 |1 - f(v)|^{-1} \cdot |f(v)|.
$$
\end{lemm}
\begin{proof}
Put $\zeta = f(v)$ and notice that $M(z) = \frac{z + \zeta}{1 + \ov{\zeta}z}$ maps $\D$ biholomorphically onto itself and is such that $M(-\zeta) = 0$ and $M(0) = \zeta$.
Applying Theorem~$1.3$ of~\cite{Pom} to $f^{-1} \circ M$ we have,
$$
|\xi - v| = |f^{-1} \circ M(- \zeta) - f^{-1} \circ M(0)|
\le |(f^{-1} \circ M)'(0)| \frac{|\zeta|}{(1 - |\zeta|)^2}.
$$
Then the lemma follows considering that $|M'(0)| = 1 - |\zeta|^2 \le 2 (1 - |\zeta|)$.
\end{proof}
\subsection{Preliminary lemmas.}\label{ss:lemmas for summ and CE}
Let $v \in \CVJ$ and $c \in \CJ$ be given.
Let $k$ be a positive integer such that $\dist(R^k(v), c) < r_K$ (where $r_K > 0$ is as in~\S\ref{ss:distortion}) and let $r > 0$ be such that $R^k(c) \in \partial \tB(c, r)$.
Then we will say that $k$ is an \textsf{\textit{univalent time for~$v$ and~$c$}}, if the pull-back of the closure of $\tB(c, r)$ to $v$ by $R^k$ is univalent.
We denote by $( k_i(v, c) )_{i \ge 1}$ the increasing sequence of all of the univalent times for~$v$ and~$c$ and for each $i \ge 1$ we denote by $\xi_i(v, c)$ the element of $R^{-k_i(v, c)}(c)$ contained in the corresponding pull-back.
\begin{lemm}\label{l:summability0}
Suppose that $R$ satisfies the summability condition with exponent $\beta \in (0, 1]$.
When $\beta = 1$ let $( \eta_j )_{j \ge 1}$ be such that $\eta_j \rightarrow \infty$ as $j \rightarrow \infty$ and such that
$$
\sum_{j \ge 1} \frac{\eta_j}{|(R^j)'(v)|^\beta} < \infty,
\text{ for every } v \in \CVJ.
$$
When $\beta < 1$ put $\eta_j = 1$, for $j \ge 1$.
Given a constant $C > 0$ put, 
$$
\rho(\delta) = C \inf_{\dist(\xi_i(v, c), v) \ge \delta} 
\left( \frac{\dist(\xi_i(v, c), v)}{\delta} \right) 
|(R^{k_i(v, c) + 1})'(v)|^{1 - \beta} \eta_{k_i(v, c) + 1}.
$$
If $C$ is small enough, then for every $c \in \CJ$, every $n \ge 1$ and every $\xi \in R^{-n}(c)$, such that $R^i(\xi) \not \in \tB(\Crit, \delta)$ for $i = 0, \ldots, n - 1$, the pull-back of $\tB(c, \delta \rho(\delta))$ to~$\xi$ by $R^n$ is univalent.
\end{lemm}
\begin{proof}
Let $D > 0$ be the constant such that the numbers
\begin{equation}\label{e:2}
d_j = D \cdot
\eta_{j + 1} \max_{v \in \CVJ}|(R^{j + 1})'(v)|^{-\beta},
\end{equation}
satisfy $\prod_{j \ge 0} (1 - d_j) = \tfrac{1}{2}$, as in the shrinking neighborhoods, and let $K > 0$ be the corresponding distortion constant as in~\eqref{e:1}.
Given $\delta > 0$ and $\dconst > 1$ such that $2\delta \dconst < r_K$, consider the shrinking neighborhoods with~$r = 2\delta \dconst$; i.e. for $j = 0, \ldots, n$ consider the pull-backs $(U_j, U_j')$ that contain $R^{n - j}(\xi)$, as explained in~\S\ref{ss:shrinking neighborhoods}.

Let $k \in \{ 0, \ldots, n \}$ be the least integer, if any, so that $U_k \cap \CV \neq \emptyset$ and let $v \in \CV$ be such that $v \in U_k$.
We assume that $v \in J(R)$, see~\S\ref{s:preliminaries}.
So $R^k$ is univalent on $U_k'$ and then there is $i \ge 1$ such that $k = k_i(v, c)$ and such that $\xi = \xi_i(v, c) \in U_k$ is the unique preimage of $c$ by $R^k$ in $U_k$.

Note that for some $C_0 > 0$ only depending on $R$ we have
$$
\frac{\dist(\xi,v)}{(\delta \dconst)^\frac{1}{\mu_{c}}} 
\le  C_0 \frac{\dist(\xi,v)}{\dist(c, R^k(v))}.
$$
Considering \eqref{e:2} and property~\eqref{e:1} of shrinking neighborhoods, this last quantity is at most
$$
C_0 \left( \frac{\dist(\xi,v)}{\dist(c, R^k(v))}
\right)^\frac{\mu_{c} - 1}{\mu_{c}} 
|(R^k)'(v)|^{-\frac{1}{\mu_{c}}} 
\left(KD^{-1}\eta_{k + 1}^{-1}
|(R^{k + 1})'(v)|^\beta 
\right)^\frac{1}{\mu_{c}}.
$$
Since $\dist(c, R^k(v))^{\mu_c - 1} \sim |R'(R^k(v))|$, for some constant $C_1 > 0$ we obtain
$$
\frac{\dist(\xi,v)}{(\delta \dconst)^\frac{1}{\mu_{c}}} 
\le C_1 \left( \dist(\xi,v) \right)^\frac{\mu_{c} - 1}{\mu_{c}} 
|(R^{k + 1})'(v)|^{\frac{\beta - 1}{\mu_{c}}}
\eta_{k + 1}^{- \frac{1}{\mu_{c}}}.
$$
Thus
$$
\delta \dconst \ge C_1^{-\mu_{c}} \dist(\xi,v) 
|(R^{k + 1})'(v)|^{1 - \beta}\eta_{k + 1}.
$$
By hypothesis $\dist(\xi, v) \ge \delta$, so letting $C = \tfrac{1}{2}C_1^{-\mu_{c}}$ in the definition of $\rho$, we have $\dconst > \rho(\delta)$.
Hence, if we take $\dconst = \rho(\delta)$, the neighborhoods $U_k$ avoid critical values and therefore the corresponding pull-back of $\tB(c, \delta \rho(\delta))$ is univalent.
\end{proof}
\begin{lemm}\label{l:summability1}
Let $\kappa_0 \in (0,1)$ be the constant as in the proof of part~$2$ of the proof of Proposition~\ref{p:BC-UP}.
Given a constant $C_0 > 0$ set
$$
\rho_1(\delta) = C_0 \inf_{\dist(\xi_i(v, c), v) < \delta} \left(\frac{\delta}{\dist(\xi_i(v, c),v)} \right)^{\mu_c - 1} |(R^{k_i(v, c) + 1})'(v)|.
$$
If $C_0$ is small enough, then the rational map $R$ is \BCg{} with function $r_0 \= \min \{ \frac{\kappa_0}{2} \rho, \rho_1 \}$.
\end{lemm}
\begin{proof}
Consider $\delta > 0$ small enough so that $r_0(\delta) \ge 2$.
We need to prove that for every $c \in \CJ$ and every sequence of pull-backs $U_0 = \tB(c, \delta r_0(\delta))$, $U_1, \ldots, U_k$ such that $U_k \cap B(\CV, \delta) \neq \emptyset$, we have $\diam(U_k) \le \delta$.
We will prove this assertion in the special case when for every $j = 0, \ldots, k - 1$ we have $U_j \cap B(\CV, \delta) = \emptyset$.
As $r_0(\delta) \ge 2$, the general case follows by induction.

So assume that for every $j = 0, \ldots, k - 1$ we have $U_j \cap B(\CV, \delta) = \emptyset$.
Consider the corresponding pull-backs $U_j'$ and $U_j''$ of $B(c, \kappa_0^{-1}\delta r_0(\delta))$ and $B(c_0, 2 \kappa_0^{-1} \delta r_0(\delta))$ respectively, so that $U_j \subset U_j' \subset U_j''$.
Since $r_0(\delta) \le \frac{\kappa_0}{2} \rho(\delta)$, we have by the previous lemma that $R^k : U_k'' \to U_0''$ is univalent.

Assume that there is $v \in U_k' \cap \CV$.
Hence there is $i \ge 0$ such that $k = k_i(v, c)$ and such that $\xi = \xi_i(v, c) \in U_k'$ is the $k$-th preimage of $c$ in $U_k'$.
Let $c_0 \in \CJ$ be so that $R(c_0) = v$ and let $U_{k + 1}'$ be the connected component of $R^{-1}(U_k')$ that contains~$c_0$.
Since $R^{k + 1}$ is not univalent on $U_{k + 1}'$, the conclusion of Lemma~\ref{l:summability0} applied to a preimage~$z$ of $\xi$ by $R$ in $U_{k + 1}'$ and with $n = k + 1$, is false.
Therefore we must have $\xi \in B(v, \delta)$.
By \KDT{} the distortion of $R^k$ on $U_k'$ is bounded by some definite constant $D > 1$, so there is $C_1 > 0$ only depending on $R$ such that,
\begin{multline*}
\frac{\delta}{(\delta r_0(\delta))^\frac{1}{\mu_c}} 
\le
C_1 \frac{\diam(U_k)}{\diam(U_0)}
\le C_1 D \frac{\dist(\xi, v)}{\dist(c, R^k(v))}
\le \\ \le
C_1 D^{1 + \frac{1}{\mu_c}}
\left(\frac{\dist(\xi, v)}{\dist(c, R^k(v))} \right)^\frac{\mu_c - 1}{\mu_c} 
|(R^k)'(v)|^\frac{1}{\mu_c}.
\end{multline*}
So there is a constant $C_2 > 0$, only depending on $R$, such that
$$
r_0(\delta)
\ge
C_2 \left(\frac{\delta}{\dist(\xi, v)}\right)^{\mu_c - 1} |(R^{k + 1})'(v)|.
$$
Since $\dist(\xi, v) < \delta$, taking $C_0$ equal to $C_2 / 2$ we obtain a contradiction.
So, $U_k' \cap \CV = \emptyset$ for this choice of $C_0$.
Therefore $\diam(U_k) \le \delta$, by definition of $\kappa_0$.
\end{proof}

\

The following lemma is needed for the Collet-Eckmann case.
\begin{lemm}\label{l:slowrecurrence}
Suppose that there is a constant $C_0 > 0$ such that for every $v \in \CVJ$ and $k \ge 1$ we have,
$$
|(R^{k + 1})'(v)| \ge C_0.
$$
There there are constants $C_1 > 0$ and $\theta \in (0, 1)$ such that for every $v \in \CVJ$, every $k \ge 1$ and every $\xi \in R^{-k}(\Crit(R))$ we have,
$$
\dist(\xi, v) \ge C_1 \theta^k.
$$
\end{lemm}
\begin{proof}
Let $M = \sup_{\CC}|R'|$ and let $C > 0$ such that for every $w \in \CC$ we have $|R'(w)| \le C \dist(z, \Crit(R))$.

From $|(R^{k + 1})'(v)| \ge C_0$ it follows that $|R'(R^k(v))| \ge C_0 M^{-k}$ so
$$
\dist(R^k(v), \Crit(R)) \ge C^{-1}C_0M^{-k}.
$$
Hence the pull-back of $B(R^k(v), C^{-1}C_0M^{-k})$ by $R^k$ to $v$ is disjoint from $R^{-k}(\Crit(R))$.
Since the later set contains $B(v, C^{-1}C_0M^{-2k})$ the lemma follows with $C_1 = C^{-1}C_0$ and $\theta = M^{-2}$.
\end{proof}
\subsection{Proof of Theorem~\ref{t:summ and CE}}\label{ss:proof of summ and CE}
\

\partn{1}
Considering that $k_i(v, c)$, and hence $\eta_{k_i(v, c)}$ and $|(R^{k_i(v, c) + 1})'(v)|$, are big when $\dist(\xi_i(v, c), v)$ is small, we conclude that the function $r_0$ defined in Lemma~\ref{l:summability1} is such that $r_0(\delta) \to \infty$ as $\delta \to 0$.

\partn{2}
Let $\theta \in (0, 1)$ be given.
By part~$3$ of Proposition~\ref{p:BC-UP} is enough to prove that the sum
$
\sum_{n \gg 1} \left( r(\theta^n) \right)^{-\alpha}
$
is finite.
But this sum is bounded from above by,
\begin{multline*}
C \sum_{i, v, c} \left( \ \sum_{\theta^n < \dist(\xi_i(v, c), v)}
\left(\frac{\dist(\xi_i(v, c), v)}{\theta^n} \right)^{- \beta / (1 - \beta)}
|(R^{k_i(v, c) + 1})'(v)|^{- \beta} \right.
\\
\left. + \sum_{\theta^n > \dist(\xi_i(v, c), v)}
\left(\frac{\theta^n}{\dist(x_i(v, c),v)} \right)^{ - (\mu_c - 1) \beta/(1 - \beta)}
|(R^{k_i(v, c) + 1})'(v)|^{ - \beta / (1 - \beta)}  \right)
\\ \le
\widetilde{C}\sum_{i, v, c} |(R^{k_i(v, c) + 1})'(v)|^{- \beta}
<
\infty.
\end{multline*}
\partn{3}
Suppose that $R$ satisfies the Collet-Eckmann condition.
That is, there are $C_0 > 0$ and $\lambda > 1$ such that for every $v \in \CVJ$ we have $|(R^k)'(v)| \ge C_0 \lambda^k$.
In particular, for every $\beta \in (0, 1)$, the rational map $R$ satisfies the summability condition with exponent $\beta$.
Fix $\beta \in (0, 1)$ and let $\rho$ and $\rho_1$ be the functions defined in lemmas~\ref{l:summability0} and~\ref{l:summability1}, respectively.

By Lemma~\ref{l:slowrecurrence} there is $C_1 > 0$ and $\theta \in (0,1)$ such that $\dist(\xi_i(v, c), v) \ge C_1 \theta^{k_i(v, c)}$.
Therefore there is $C_2 > 0$ and $\gamma \in (0,1)$ such that,
$$
|(R^{k_i(c, v) + 1})'(v)|
\ge
C_0 \lambda^{k_i(v, c)}
\ge
C_2(\dist(\xi_i(v, c), v))^{-\gamma}.
$$ 
Choose $\mu \in (0,1)$.
If $\dist(\xi_i(v, c), v) \le \delta^{1 - \mu}$, we have
$$
|(R^{k_i(v, c) + 1})'(v)| \ge
C_2 \dist(\xi_i(v, c), v)^{-\gamma}
\ge C_2 \delta^{-\gamma(1 - \mu)};
$$
in particular $\rho_1(\delta) \ge C_0C_2 \delta^{-\gamma(1 - \mu)}$.

On the other hand note that
$$
\text{if } \
\dist(\xi_i(v, c), v) > \delta^{1 - \mu}
\ \text{ then } \
\dist(\xi_i(v, c), v) / \delta > \delta^{-\mu}.
$$
Therefore $\rho(\delta) \ge C \min \{ C_0\delta^{-\mu}, C_2 \delta^{-\gamma(1 - \mu)(1 - \beta)} \}$ in this case.
Thus 
$$
r_0(\delta)
=
\min\{ \frac{\kappa_0}{2}\rho(\delta), \rho_1(\delta) \}
\ge
C_4 \delta^{-\alpha},
$$
where $\alpha = \min \{ \mu, \gamma(1 - \mu)(1 - \beta)\}$ and $C_4$ is some positive constant.

\partn{4}
If the closure of the forward orbit of each critical value does not contain critical points, then clearly there is $C > 0$ such that $R$ satisfies the Univalent Pull-back Condition with function $r(\delta) \= C \delta^{-1}$.
So, Part~$4$ of Theorem~\ref{t:summ and CE} follows from Proposition~\ref{p:BC-UP}.

\section{Quasi-conformal homeomorphisms.}\label{a:qc}
In this appendix we review some properties of quasi-conformal maps.
See~\cite{LV} and~\cite{Ah} for references.
\subsection{Modulus of annuli.}
Every annulus $A \subset \CC$ is either conformally equivalent to $\C^* \= \C \setminus \{ 0 \}$ or to $\{ z \in \C \mid 1 < |z| < R \}$, for some constant $R \in (1, + \infty]$.
In the latter case, the constant $R$ is uniquely determined by $A$ and then $\modulus(A) \= \ln(R)$ is called the \textsf{\textit{modulus}} of $A$.
(Some authors call $\frac{1}{2\pi} \ln R$ the modulus of $A$; here we follow~\cite{LV}.)
The modulus $\modulus(A)$ can also be defined by,
\begin{equation}\label{e:def modulus}
\modulus(A)
=
\left( \frac{1}{2\pi} \inf_{h} \iint_A |\nabla h |^2 dxdy \right)^{-1},
\end{equation}
where the infimum is taken over all functions $h : A \rightarrow (0, 1)$ of class $C^1$, such that $h(z) \rightarrow 1$ as $z$ approaches an end of $A$ and $h(z) \rightarrow 0$ as $z$ approaches the other end of $A$.

Consider the \textsf{\textit{flat metric}} $|\frac{dz}{z}|$ on $\C^*$, that makes $\C^*$ isometric to $S^1 \times \R$.
The \textsf{\textit{flat metric}} on an arbitrary annulus is defined by pulling back the flat metric on $\C^*$, by a conformal representation of the annulus to $\C^*$ or to $\{ z \in \C \mid 1 < |z| < R \}$.
\subsection{Quasi-conformal homeomorphisms.}
Given a constant $K \ge 1$, we say that a homeomorphism $\chi$ between open subsets of $\CC$ is $K$ \textsf{\textit{quasi-conformal}} or $K$-\textsf{\textit{qc}}, if the following equivalent conditions hold.
\begin{enumerate}
\item[1.]
For every annulus $A$ contained in the domain of $\chi$ we have
$$
K^{-1} \modulus(A) \le \modulus(\chi(A)) \le K \modulus(A).
$$
\item[2.]
The homeomorphism $\chi$ has a distributional derivative that is locally $L^2$ and $\| D \chi \|^2 \le K \Jac(D\chi)$, on a set of full Lebesgue measure.
\end{enumerate} 
The constant $K \ge 1$ is called the \textsf{\textit{dilatation}} of $\chi$.
If we do not want to specify the dilatation, we just say that $\chi$ is \textsf{\textit{quasi-conformal}} or \textsf{\textit{qc}}.
Note that every conformal homeomorphism is $1$-qc.
Conversely, every $1$-qc homeomorphism is conformal.

Property~$1$ above implies that the inverse of a $K$-qc homeomorphism is also a $K$-qc homeomorphism.
For every qc homeomorphism there is a set of full Lebesgue measure where it is differentiable (in the usual sense) and where its derivative coincides with its distributional derivative.


The following lemma can be found in~\cite{DH}.
\begin{generic}[Gluing Lemma]
Let $U$ be a bounded open subset of $\C$ and let $\chi : \CC \to \CC$ be a qc homeomorphism.
Let $\chi_0 : U \to \chi(U)$ be a qc homeomorphism such that the map $\tchi : \CC \to \CC$ that is equal to $\chi$ outside $U$ and is equal to $\chi_0$ on $U$ is continuous.
Then the following properties hold.
\begin{enumerate}
\item[1.]
$\tchi$ is a qc homeomorphism of $\CC$.
\item[2.]
The derivatives of $\tchi$ and $\chi$ coincide on a set of full Lebesgue in $\CC \setminus U$.
\end{enumerate}
\end{generic}
Qc homeomorphisms preserve sets of Lebesgue measure zero and sets of $\sigma$-finite length qc removable: if $\chi: U \to \chi(U)$ is a homeomorphism that is $K$-qc outside a set of $\sigma$-finite length, then $\chi$ is $K$-qc.
\subsection{Normalized qc maps.}
Let $\cF$ be the collection of homeomorphisms of $\CC$ or the collection of homeomorphisms of $\D$.
We say that an element of $\cF$ is \textsf{\textit{normalized}} if it fixes~3 base points in $\CC$ (resp, if it fixes $0$).
In each case, for each element $h$ of $\cF$ there exist a M\"obius transformation $\varphi$ (fixing $\D$ in the latter case) such that $\varphi \circ h$ is a normalized element of $\cF$.

For each constant $K \ge 1$, the collection of normalized $K$-qc homeomorphisms in~$\cF$ is compact: Every sequence of normalized $K$-qc homeomorphisms in~$\cF$ admits a subsequence that converges uniformly to a normalized $K$-qc homeomorphism in~$\cF$.

\begin{lemm}\label{l:approximation conformal}
Let $K \ge 1$ be given and let $\cF$ be one of the collections described above.
Then the following assertions hold.
\begin{enumerate}
\item[1.]
Let $( \varepsilon_k )_{k \ge 1}$ be a sequence of positive numbers converging to~0 and let $( \chi_k )_{k \ge 1}$ be a sequence of normalized $K$-qc homeomorphisms in $\cF$ that converges uniformly to a homeomorphism $\chi$ and such that $\chi_k$ is conformal outside a set of Lebesgue measure $\varepsilon_k$.
Then $\chi$ is conformal.
\item[2.]
For every $\varepsilon > 0$ there is $\delta > 0$ such that every normalized $K$-qc homeomorphism of $\CC$ in $\cF$ that is conformal outside a set of Lebesgue measure $\delta$, is~$\varepsilon$ close to a normalized conformal map.
\end{enumerate}
\end{lemm}

The proof of this lemma is based on the following one.
\begin{lemm}\label{l:mod}
Let $A$ be an annulus and let $\chi : A \to \chi(A)$ be a $K$-qc homeomorphism.
Let ${\mathcal N}$ be a subset of $A$ such that $\chi$ is conformal on $A \setminus {\mathcal N}$. 
Then,
$$
\modulus(A) \left(1 + \frac{(K - 1)|{\mathcal N}|}{2\pi \modulus(A)} \right)^{-1}
\le \modulus(\chi(A)) \le \modulus(A) + \frac{(K - 1)|{\mathcal N}|}{2\pi},
$$
where $|{\mathcal N}|$ denotes the area of ${\mathcal N}$ with respect to the flat metric of $A$.
\end{lemm}
\begin{proof}
For the upper bound see~\cite{LV}~($6.6$), p.~$221$.
We will now prove the lower bound.

Suppose that $A$ is the straight cylinder $S^1 \times (0, \modulus(A))$ and denote by $\pi$ the projection $\pi: A \to (0,1)$ given by $\pi(\theta, t) = \frac{t}{\modulus(A)}$.
By~\eqref{e:def modulus} we have,
\begin{multline*}
\frac{1}{\modulus(\chi(A))}
\le
\frac{1}{2\pi} \iint_{\chi(A)} |\nabla(\pi \circ \chi^{-1})|^2 \ dxdy
\\ \le
\frac{1}{2\pi}(\modulus(A))^{-2} \iint_{\chi(A)} \| D\chi^{-1} \|^2 \ dxdy
\\ \le
\frac{1}{2\pi}(\modulus(A))^{-2} \left( \iint_{\chi(A)} \Jac(\chi^{-1}) \ dxdy
+ (K - 1)\iint_{\chi({\mathcal N})}  \Jac(\chi^{-1}) \ dx dy \right)
\\
=
\frac{2\pi \modulus(A) + (K - 1)|{\mathcal N}|}{2\pi (\modulus(A))^2}.
\end{multline*}
\end{proof}
\begin{proof}[Proof of Lemma~\ref{l:approximation conformal}.]
\

\partn{1}
For a given annulus $A$, the previous lemma implies that $\modulus(\chi_k(A)) \rightarrow \modulus(A)$ as $k \to \infty$.
We conclude that $\modulus(\chi(A)) = \modulus(A)$ and that $\chi$ preserves the modulus of every annulus.
So $\chi$ is $1$-qc, and hence conformal.

\partn{2}
If this is not true, then there is $\varepsilon > 0$ and a sequence of normalized $K$-qc homeomorphisms $( \chi_k )_{k \ge 1}$ such that $\chi_k$ is conformal outside a set of measure $\frac{1}{k}$ and such that the distance from $\chi_k$ to each normalized conformal map is at least~$\varepsilon$.
This contradicts part~$1$.
\end{proof}

\bibliographystyle{plain}

\begin{thebibliography}{GSm2}
\bibitem[Ah]{Ah} L. V. Ahlfors. {\it Lectures on quasiconformal mappings.}  
Van Nostrand, 1966.

\bibitem[Asp]{Asp} M. Aspenberg.
{\it The Collet-Eckmann condition for rational functions on the Riemann sphere.}
Doctoral Thesis, KTH, Sweden 2004.

\bibitem[Ast]{As} K. Astala. {\it Area distortion property of quasi-conformal maps.}
Acta Math. {\bf 173}  (1994),  37-60.

\bibitem[Av]{Av} A. Avila.
{\it Infinitesimal perturbations of rational maps.}
Nonlinearity {\bf 15} (2002), 695-704.

\bibitem[BS]{BvS} H. Bruin, S. van Strien.
{\it Expansion of derivatives in one dimentional dynamics.}
Israel J. Math. {\bf 137}  (2003), 223-263.


\bibitem[CG]{CG} L. Carleson, T. Gamelin. {\it Complex dynamics}.
Springer-Verlag, 1993.

\bibitem[CE]{CE83} P. Collet, J.P. Eckmann. {\it Positive Liapunov exponents and absolute continuity for maps of the interval.} Ergodic Theory and Dynamical Sys. {\bf 3} (1983), 13-46.

\bibitem[DH1]{DH} A. Douady and J. Hubbard.
{\it On the dynamics of polynomial-like mappings}.
Ann. Sci. Ec. Norm. Sup. {\bf 18} (1985), 287-344.

\bibitem[DH2]{DHalgorithm}
A. Douady, J.H. Hubbard.
{\it A proof of Thurston's topological characterization of rational functions.}
Acta Math. {\bf 171} (1993),  263-297.

\bibitem[DU]{DU} M. Denker, M. Urba\'nski.
{\it Hausdorff and conformal measures on Julia sets with rationally indifferent periodic points.}
J. London Math. Soc. {\bf 43} (1991), 107-118.

\bibitem[GSm1]{GSmCE} J. Graczyk, S. Smirnov. {\it Collet, Eckmann and Holder}.
Invet. math. {\bf 133.}, (1998), 69-96.

\bibitem[GSm2]{GSm} J. Graczyk, S. Smirnov. {\it Weak expansion and geometry of Julia sets}. March 1999 version.

\bibitem[GSw]{GSw} J. Graczyk, G. \'Swi\c atek.
{\it Harmonic measure and expansion on the boundary of the connectedness locus.}
Invent. Math.  {\bf 142}  (2000), 605-629.

\bibitem[Ha{\"i}]{Hai} P. Ha{\"i}ssinsky.
{\it Pincement de polyn\^omes.}
Comment. Math. Helv. {\bf 77} (2002), 1-23.

\bibitem[HS]{HS} J.H. Hubbard, D. Schleicher.
{\it The spider algorithm.}
Complex dynamical systems (Cincinnati, OH, 1994),  155-180, Proc. Sympos. Appl. Math., 49, Amer. Math. Soc., Providence, RI, 1994.

\bibitem[K]{Khan} J. Khan.
{\it Holomorphic removability of Julia sets.}
IMS preprint at Stony Brook \#1998/11.

\bibitem[LV]{LV} O. Lehto, K.I. Virtanen. {\it Quasiconformal mappings in the plane}.
Springer-Verlag 1973.

\bibitem[L]{Levsumm} G. Levin.
{\it On an analytic approach to the Fatou conjecture.}
Fund. Math. {\bf 171} (2002), 177-196.

\bibitem[Mak]{Mak} P. Makienko.
{\it Remarks on the Ruelle operator and the invariant line fields problem. II.}
Ergodic Theory Dynam. Systems  {\bf 25} (2005), 1561-1581.

\bibitem[Ma\~n]{Mane} R. Ma\~n\'e.
{\it On a theorem of Fatou.}
Bol. Soc. Brasil. Mat. {\bf 24} (1993), 1-11.

\bibitem[MR]{MR} R. Ma\~n\'e, L.F. da Rocha.
{\it Julia sets are uniformly perfect.}
 Proc. Amer. Math. Soc.  {\bf 116}  (1992), 251-257.

\bibitem[MSS]{MSS} R. Ma\~n\'e, P. Sad, D. Sullivan.
{\it On the dynamics of rational maps}.
Ann. Sci. Ec. Norm. Sup. {\bf 16} (1983), 193-217.

\bibitem[Mar]{Mar} M. Martens. {\it Distortion results and invariant Cantor sets of unimodal mappings}. Erg. Th. and Dyn. Sys. {\bf 14} (1994), 331-349.

\bibitem[McM]{McM} C. McMullen. {\it Renormalization and $3$-manifolds
which fiber over the circle}.
Annals of mathematics studies, {\bf 142}, Princeton University Press, 1996.

\bibitem[Mil]{Milbook} J. Milnor. {\it Dynamics in one complex variable. Introductory lectures.}
Friedr. Vieweg \& Sohn, Braunschweig, 1999.


\bibitem[NS]{NS} T. Nowicki, S. van Strein. {\it Absolutely continuous invariant measures under the summability condition.}
Invent. Math. {\bf 105} (1991), 123-136.

\bibitem[Pil]{Pil} K. Pilgrim.
{\it Bounded geometry of quadrilaterals and variation of multipliers for rational maps.}
Fund. Math. {\bf 182} (2004), 137-150.

\bibitem[Po]{Pom} C. Pommerenke. {\it Boundary behavior of conformal maps.} Springer-Verlag, Berlin, 1992.

\bibitem[Pr1]{Pr2} F. Pryzytycki.
{\it On measure and Hausdorff dimension of Julia sets for holomorphic Collet-Eckmann maps}.
Int. conf. on dynamical systems, Montevideo 1995- a tribute to R. Ma\~n\'e,
Eds. F. Ledrappier, J. Lewowicz, S. Newhouse.
Pitman Res. Notes in Math. series, 362. Longman 1996, 167-181.

\bibitem[Pr2]{Pr1} F. Przytycki. {\it Iterations of holomorphic Collet-Eckmann maps:
Conformal and Invariant measures.
Appendix: On non-renormalizable quadratic polynomials}.
Trans. A.M.S. {\bf 350} (1998), 717-742.


\bibitem[PR]{PR} F. Przytycki, J. Rivera-Letelier.
Statistical properties of Topological Collet-Eckmann maps.
Preprint February~2006.

\bibitem[PRS]{PRS} F. Przytycki, J. Rivera-Letelier, S. Smirnov.
{\it Equivalence and topological invariance of conditions for non-uniform hyperbolicity in iteration of rational maps.}
Invent. Math.  {\bf 151}  (2003), 29-63.

\bibitem[PU]{PrUparabolic} F. Przytycki, M. Urba\'nski.
{\it Porosity of Julia sets of non-recurrent and parabolic Collet Eckmann functions.}
Ann. Acad. Fenn. {\bf 26} (2001), 125-154.

\bibitem[R]{ims} J. Rivera-Letelier. {\it Rational maps with decay of geometry: rigidity, Thurston's algorithm and local connectivity.} IMS preprint at Stony Brook \#2000/11.

\bibitem[RS]{RS} J. Rivera-Letelier, S. van Strien.
{\it Decay of geometry and growth of derivatives.}
September 2001.

\bibitem[Se]{Sester}
O. Sester.
{\it Dynamique des Polyn\^omes Fibr\'es.}
Thesis, Universit\'e Paris-Sud~1997.

\bibitem[SL]{SL} M. Shishikura, Tan Lei.
{\it An alternative proof of Ma\~n\'e's theorem on non-expanding Julia sets.}  The Mandelbrot set, theme and variations,  265-279, London Math. Soc. Lecture Note Ser., 274, Cambridge Univ. Press, Cambridge, 2000.

\bibitem[Sm]{SmCE} S. Smirnov. {\it Symbolic dynamics and Collet-Eckmann condition}.
Internat. Math. Res. Notices {\bf 7} (2000), 333-351.

\bibitem[Sl]{Sl} Z. Slodkowsky. {\it Holomorphic motions and polynomials hulls}. 
Proc. Amer. Math. Soc. {\bf 111} (1991), 347-355.

\bibitem[U]{U} M. Urba\'nski.
{\it Measures and dimensions in conformal dynamics.}
Bull. Amer. Math. Soc. {\bf 40}  (2003), 281-321

\end{thebibliography}

\end{document}